\ifx\newheadisloaded\relax\immediate\write16{***already loaded}\endinput\else\let\newheadisloaded=\relax\fi
\gdef\islinuxolivetti{F}
\gdef\PSfonts{T}
\magnification\magstep1

\newdimen\papwidth
\newdimen\papheight
\newskip\beforesectionskipamount  
\newskip\sectionskipamount 
\def\sectionskip{\vskip\sectionskipamount}
\def\beforesectionskip{\vskip\beforesectionskipamount}
\papwidth=16truecm
\if F\islinuxolivetti
\papheight=22truecm
\voffset=0.4truecm
\hoffset=0.4truecm
\else
\papheight=22truecm
\voffset=-1.5truecm
\hoffset=-2.3truecm
\fi
\hsize=\papwidth
\vsize=\papheight
\catcode`\@=11
\ifx\amstexloaded@\relax
\else
\nopagenumbers
\headline={\ifnum\pageno>1 {\hss\tenrm-\ \folio\ -\hss} \else
{\hfill}\fi}
\fi
\catcode`\@=\active
\newdimen\texpscorrection
\texpscorrection=0.15truecm 

\def\sectionsize{\twelvepoint}
\def\sectiontype{\bf}
\def\subsectionsize{}
\def\subsectiontype{\bf}
\def\em{\sl}
\newfam\truecmsy
\newfam\truecmr
\newfam\msbfam
\newfam\scriptfam
\newfam\frakfam
\newfam\frakbfam

\newskip\ttglue 
\if T\islinuxolivetti
\papheight=11.5truecm
\fi
\if F\PSfonts
\font\twelverm=cmr12
\font\tenrm=cmr10
\font\eightrm=cmr8
\font\sevenrm=cmr7
\font\sixrm=cmr6
\font\fiverm=cmr5

\font\twelvebf=cmbx12
\font\tenbf=cmbx10
\font\eightbf=cmbx8
\font\sevenbf=cmbx7
\font\sixbf=cmbx6
\font\fivebf=cmbx5

\font\twelveit=cmti12
\font\tenit=cmti10
\font\eightit=cmti8
\font\sevenit=cmti7
\font\sixit=cmti6
\font\fiveit=cmti5

\font\twelvesl=cmsl12
\font\tensl=cmsl10
\font\eightsl=cmsl8
\font\sevensl=cmsl7
\font\sixsl=cmsl6
\font\fivesl=cmsl5

\font\twelvei=cmmi12
\font\teni=cmmi10
\font\eighti=cmmi8
\font\seveni=cmmi7
\font\sixi=cmmi6
\font\fivei=cmmi5

\font\twelvesy=cmsy10	at	12pt
\font\tensy=cmsy10
\font\eightsy=cmsy8
\font\sevensy=cmsy7
\font\sixsy=cmsy6
\font\fivesy=cmsy5
\font\twelvetruecmsy=cmsy10	at	12pt
\font\tentruecmsy=cmsy10
\font\eighttruecmsy=cmsy8
\font\seventruecmsy=cmsy7
\font\sixtruecmsy=cmsy6
\font\fivetruecmsy=cmsy5

\font\twelvetruecmr=cmr12
\font\tentruecmr=cmr10
\font\eighttruecmr=cmr8
\font\seventruecmr=cmr7
\font\sixtruecmr=cmr6
\font\fivetruecmr=cmr5

\font\twelvebf=cmbx12
\font\tenbf=cmbx10
\font\eightbf=cmbx8
\font\sevenbf=cmbx7
\font\sixbf=cmbx6
\font\fivebf=cmbx5

\font\twelvett=cmtt12
\font\tentt=cmtt10
\font\eighttt=cmtt8

\font\twelveex=cmex10	at	12pt
\font\tenex=cmex10

\font\twelvemsb=msbm10	at	12pt
\font\tenmsb=msbm10
\font\eightmsb=msbm8
\font\sevenmsb=msbm7
\font\sixmsb=msbm6
\font\fivemsb=msbm5


\font\tenfrm=eufm10
\font\eightfrm=eufm8
\font\sevenfrm=eufm7
\font\sixfrm=eufm6
\font\fivefrm=eufm5

\font\tenfrb=eufb10
\font\eightfrb=eufb8
\font\sevenfrb=eufb7
\font\sixfrb=eufb6
\font\fivefrb=eufb5
\font\twelvescr=eusm10 at 12pt
\font\tenscr=eusm10
\font\eightscr=eusm8
\font\sevenscr=eusm7
\font\sixscr=eusm6
\font\fivescr=eusm5
\fi
\if T\PSfonts
\font\twelverm=ptmr	at	12pt
\font\tenrm=ptmr	at	10pt
\font\eightrm=ptmr	at	8pt
\font\sevenrm=ptmr	at	7pt
\font\sixrm=ptmr	at	6pt
\font\fiverm=ptmr	at	5pt

\font\twelvebf=ptmb	at	12pt
\font\tenbf=ptmb	at	10pt
\font\eightbf=ptmb	at	8pt
\font\sevenbf=ptmb	at	7pt
\font\sixbf=ptmb	at	6pt
\font\fivebf=ptmb	at	5pt

\font\twelveit=ptmri	at	12pt
\font\tenit=ptmri	at	10pt
\font\eightit=ptmri	at	8pt
\font\sevenit=ptmri	at	7pt
\font\sixit=ptmri	at	6pt
\font\fiveit=ptmri	at	5pt

\font\twelvesl=ptmro	at	12pt
\font\tensl=ptmro	at	10pt
\font\eightsl=ptmro	at	8pt
\font\sevensl=ptmro	at	7pt
\font\sixsl=ptmro	at	6pt
\font\fivesl=ptmro	at	5pt

\font\twelvei=cmmi12
\font\teni=cmmi10
\font\eighti=cmmi8
\font\seveni=cmmi7
\font\sixi=cmmi6
\font\fivei=cmmi5

\font\twelvesy=cmsy10	at	12pt
\font\tensy=cmsy10
\font\eightsy=cmsy8
\font\sevensy=cmsy7
\font\sixsy=cmsy6
\font\fivesy=cmsy5
\font\twelvetruecmsy=cmsy10	at	12pt
\font\tentruecmsy=cmsy10
\font\eighttruecmsy=cmsy8
\font\seventruecmsy=cmsy7
\font\sixtruecmsy=cmsy6
\font\fivetruecmsy=cmsy5

\font\twelvetruecmr=cmr12
\font\tentruecmr=cmr10
\font\eighttruecmr=cmr8
\font\seventruecmr=cmr7
\font\sixtruecmr=cmr6
\font\fivetruecmr=cmr5

\font\twelvebf=cmbx12
\font\tenbf=cmbx10
\font\eightbf=cmbx8
\font\sevenbf=cmbx7
\font\sixbf=cmbx6
\font\fivebf=cmbx5

\font\twelvett=cmtt12
\font\tentt=cmtt10
\font\eighttt=cmtt8

\font\twelveex=cmex10	at	12pt
\font\tenex=cmex10

\font\twelvemsb=msbm10	at	12pt
\font\tenmsb=msbm10
\font\eightmsb=msbm8
\font\sevenmsb=msbm7
\font\sixmsb=msbm6
\font\fivemsb=msbm5


\font\tenfrm=eufm10
\font\eightfrm=eufm8
\font\sevenfrm=eufm7
\font\sixfrm=eufm6
\font\fivefrm=eufm5

\font\tenfrb=eufb10
\font\eightfrb=eufb8
\font\sevenfrb=eufb7
\font\sixfrb=eufb6
\font\fivefrb=eufb5
\font\twelvescr=eusm10 at 12pt
\font\tenscr=eusm10
\font\eightscr=eusm8
\font\sevenscr=eusm7
\font\sixscr=eusm6
\font\fivescr=eusm5
\fi
\def\eightpoint{\def\rm{\fam0\eightrm}%
\textfont0=\eightrm
  \scriptfont0=\sixrm
  \scriptscriptfont0=\fiverm 
\textfont1=\eighti
  \scriptfont1=\sixi
  \scriptscriptfont1=\fivei 
\textfont2=\eightsy
  \scriptfont2=\sixsy
  \scriptscriptfont2=\fivesy 
\textfont3=\tenex
  \scriptfont3=\tenex
  \scriptscriptfont3=\tenex 
\textfont\itfam=\eightit
  \scriptfont\itfam=\sixit
  \scriptscriptfont\itfam=\fiveit 
  \def\it{\fam\itfam\eightit}%
\textfont\slfam=\eightsl
  \scriptfont\slfam=\sixsl
  \scriptscriptfont\slfam=\fivesl 
  \def\sl{\fam\slfam\eightsl}%
\textfont\ttfam=\eighttt
  \def\tt{\fam\ttfam\eighttt}%
\textfont\bffam=\eightbf
  \scriptfont\bffam=\sixbf
  \scriptscriptfont\bffam=\fivebf
  \def\bf{\fam\bffam\eightbf}%
\textfont\frakfam=\eightfrm
  \scriptfont\frakfam=\sixfrm
  \scriptscriptfont\frakfam=\fivefrm
  \def\frak{\fam\frakfam\eightfrm}%
\textfont\frakbfam=\eightfrb
  \scriptfont\frakbfam=\sixfrb
  \scriptscriptfont\frakbfam=\fivefrb
  \def\bfrak{\fam\frakbfam\eightfrb}%
\textfont\scriptfam=\eightscr
  \scriptfont\scriptfam=\sixscr
  \scriptscriptfont\scriptfam=\fivescr
  \def\script{\fam\scriptfam\eightscr}%
\textfont\msbfam=\eightmsb
  \scriptfont\msbfam=\sixmsb
  \scriptscriptfont\msbfam=\fivemsb
  \def\bb{\fam\msbfam\eightmsb}%
\textfont\truecmr=\eighttruecmr
  \scriptfont\truecmr=\sixtruecmr
  \scriptscriptfont\truecmr=\fivetruecmr
  \def\truerm{\fam\truecmr\eighttruecmr}%
\textfont\truecmsy=\eighttruecmsy
  \scriptfont\truecmsy=\sixtruecmsy
  \scriptscriptfont\truecmsy=\fivetruecmsy
\tt \ttglue=.5em plus.25em minus.15em 
\normalbaselineskip=9pt
\setbox\strutbox=\hbox{\vrule height7pt depth2pt width0pt}%
\normalbaselines
\rm
}

\def\tenpoint{\def\rm{\fam0\tenrm}%
\textfont0=\tenrm
  \scriptfont0=\sevenrm
  \scriptscriptfont0=\fiverm 
\textfont1=\teni
  \scriptfont1=\seveni
  \scriptscriptfont1=\fivei 
\textfont2=\tensy
  \scriptfont2=\sevensy
  \scriptscriptfont2=\fivesy 
\textfont3=\tenex
  \scriptfont3=\tenex
  \scriptscriptfont3=\tenex 
\textfont\itfam=\tenit
  \scriptfont\itfam=\sevenit
  \scriptscriptfont\itfam=\fiveit 
  \def\it{\fam\itfam\tenit}%
\textfont\slfam=\tensl
  \scriptfont\slfam=\sevensl
  \scriptscriptfont\slfam=\fivesl 
  \def\sl{\fam\slfam\tensl}%
\textfont\ttfam=\tentt
  \def\tt{\fam\ttfam\tentt}%
\textfont\bffam=\tenbf
  \scriptfont\bffam=\sevenbf
  \scriptscriptfont\bffam=\fivebf
  \def\bf{\fam\bffam\tenbf}%
\textfont\frakfam=\tenfrm
  \scriptfont\frakfam=\sevenfrm
  \scriptscriptfont\frakfam=\fivefrm
  \def\frak{\fam\frakfam\tenfrm}%
\textfont\frakbfam=\tenfrb
  \scriptfont\frakbfam=\sevenfrb
  \scriptscriptfont\frakbfam=\fivefrb
  \def\bfrak{\fam\frakbfam\tenfrb}%
\textfont\scriptfam=\tenscr
  \scriptfont\scriptfam=\sevenscr
  \scriptscriptfont\scriptfam=\fivescr
  \def\script{\fam\scriptfam\tenscr}%
\textfont\msbfam=\tenmsb
  \scriptfont\msbfam=\sevenmsb
  \scriptscriptfont\msbfam=\fivemsb
  \def\bb{\fam\msbfam\tenmsb}%
\textfont\truecmr=\tentruecmr
  \scriptfont\truecmr=\seventruecmr
  \scriptscriptfont\truecmr=\fivetruecmr
  \def\truerm{\fam\truecmr\tentruecmr}%
\textfont\truecmsy=\tentruecmsy
  \scriptfont\truecmsy=\seventruecmsy
  \scriptscriptfont\truecmsy=\fivetruecmsy
\tt \ttglue=.5em plus.25em minus.15em 
\normalbaselineskip=12pt
\setbox\strutbox=\hbox{\vrule height8.5pt depth3.5pt width0pt}%
\normalbaselines
\rm
}

\def\twelvepoint{\def\rm{\fam0\twelverm}%
\textfont0=\twelverm
  \scriptfont0=\tenrm
  \scriptscriptfont0=\eightrm 
\textfont1=\twelvei
  \scriptfont1=\teni
  \scriptscriptfont1=\eighti 
\textfont2=\twelvesy
  \scriptfont2=\tensy
  \scriptscriptfont2=\eightsy 
\textfont3=\twelveex
  \scriptfont3=\twelveex
  \scriptscriptfont3=\twelveex 
\textfont\itfam=\twelveit
  \scriptfont\itfam=\tenit
  \scriptscriptfont\itfam=\eightit 
  \def\it{\fam\itfam\twelveit}%
\textfont\slfam=\twelvesl
  \scriptfont\slfam=\tensl
  \scriptscriptfont\slfam=\eightsl 
  \def\sl{\fam\slfam\twelvesl}%
\textfont\ttfam=\twelvett
  \def\tt{\fam\ttfam\twelvett}%
\textfont\bffam=\twelvebf
  \scriptfont\bffam=\tenbf
  \scriptscriptfont\bffam=\eightbf
  \def\bf{\fam\bffam\twelvebf}%
\textfont\scriptfam=\twelvescr
  \scriptfont\scriptfam=\tenscr
  \scriptscriptfont\scriptfam=\eightscr
  \def\script{\fam\scriptfam\twelvescr}%
\textfont\msbfam=\twelvemsb
  \scriptfont\msbfam=\tenmsb
  \scriptscriptfont\msbfam=\eightmsb
  \def\bb{\fam\msbfam\twelvemsb}%
\textfont\truecmr=\twelvetruecmr
  \scriptfont\truecmr=\tentruecmr
  \scriptscriptfont\truecmr=\eighttruecmr
  \def\truerm{\fam\truecmr\twelvetruecmr}%
\textfont\truecmsy=\twelvetruecmsy
  \scriptfont\truecmsy=\tentruecmsy
  \scriptscriptfont\truecmsy=\eighttruecmsy
\tt \ttglue=.5em plus.25em minus.15em 
\setbox\strutbox=\hbox{\vrule height7pt depth2pt width0pt}%
\normalbaselineskip=15pt
\normalbaselines
\rm
}
%
\fontdimen16\tensy=2.7pt
\fontdimen13\tensy=4.3pt
\fontdimen17\tensy=2.7pt
\fontdimen14\tensy=4.3pt
\fontdimen18\tensy=4.3pt
\fontdimen16\eightsy=2.7pt
\fontdimen13\eightsy=4.3pt
\fontdimen17\eightsy=2.7pt
\fontdimen14\eightsy=4.3pt
\fontdimen18\sevensy=4.3pt
\fontdimen16\sevensy=1.8pt
\fontdimen13\sevensy=4.3pt
\fontdimen17\sevensy=2.7pt
\fontdimen14\sevensy=4.3pt
\fontdimen18\sevensy=4.3pt
%
\def\hexnumber#1{\ifcase#1 0\or1\or2\or3\or4\or5\or6\or7\or8\or9\or
 A\or B\or C\or D\or E\or F\fi}
\mathcode`\=="3\hexnumber\truecmr3D
\mathchardef\not="3\hexnumber\truecmsy36
\mathcode`\+="2\hexnumber\truecmr2B
\mathcode`\(="4\hexnumber\truecmr28
\mathcode`\)="5\hexnumber\truecmr29
\mathcode`\!="5\hexnumber\truecmr21
\mathcode`\(="4\hexnumber\truecmr28
\mathcode`\)="5\hexnumber\truecmr29

\def\tilde{\mathaccent"0\hexnumber\truecmr7E }

\def\hat{\mathaccent"0\hexnumber\truecmr5E }
\def\dot{\mathaccent"0\hexnumber\truecmr5F }
\def\Phi{\mathchar"0\hexnumber\truecmr08 }
\def\Gamma {\mathchar"0\hexnumber\truecmr00 }
\def\Delta {\mathchar"0\hexnumber\truecmr01 }
\def\Theta {\mathchar"0\hexnumber\truecmr02 }
\def\Lambda{\mathchar"0\hexnumber\truecmr03 }
\def\Xi {\mathchar"0\hexnumber\truecmr04 }
\def\Pi{\mathchar"0\hexnumber\truecmr05 }
\def\Sigma{\mathchar"0\hexnumber\truecmr06 }
\def\Upsilon {\mathchar"0\hexnumber\truecmr07 }
\def\Phi {\mathchar"0\hexnumber\truecmr08 }
\def\Psi {\mathchar"0\hexnumber\truecmr09 }
\def\Omega{\mathchar"0\hexnumber\truecmr0A }
\newcount\EQNcount \EQNcount=1
\newcount\CLAIMcount \CLAIMcount=1
\newcount\SECTIONcount \SECTIONcount=0
\newcount\SUBSECTIONcount \SUBSECTIONcount=1
\def\ifff#1#2#3{\ifundefined{#1#2}%
\expandafter\xdef\csname #1#2\endcsname{#3}\else%
\fi}
\def\NEWDEF#1#2#3{\ifff{#1}{#2}{#3}}
\def\actualnumber{\number\SECTIONcount}
\def\EQ#1{\lmargin{#1}\eqno\tageck{#1}}
\def\NR#1{&\lmargin{#1}\tageck{#1}\cr}  
\def\tageck#1{\lmargin{#1}({\rm \actualnumber}.\number\EQNcount)
 \NEWDEF{e}{#1}{(\actualnumber.\number\EQNcount)}
\global\advance\EQNcount by 1
}
\def\SECT#1#2{\lmargin{#1}\SECTION{#2}%
\NEWDEF{s}{#1}{\actualnumber}
}
\def\SUBSECT#1#2{\lmargin{#1}
\SUBSECTION{#2} 
\NEWDEF{s}{#1}{\actualnumber.\number\SUBSECTIONcount}
}
\def\CLAIM#1#2#3\par{
\vskip.1in\medbreak\noindent
{\lmargin{#2}\bf #1\ \actualnumber.\number\CLAIMcount.} {\sl #3}\par
\NEWDEF{c}{#2}{#1\ \actualnumber.\number\CLAIMcount}
\global\advance\CLAIMcount by 1
\ifdim\lastskip<\medskipamount
\removelastskip\penalty55\medskip\fi}
\def\CLAIMNONR #1#2#3\par{
\vskip.1in\medbreak\noindent
{\lmargin{#2}\bf #1.} {\sl #3}\par
\NEWDEF{c}{#2}{#1}
\global\advance\CLAIMcount by 1
\ifdim\lastskip<\medskipamount
\removelastskip\penalty55\medskip\fi}
\def\SECTION#1{\vskip0pt plus.3\vsize\penalty-75
    \vskip0pt plus -.3\vsize
    \global\advance\SECTIONcount by 1
    \beforesectionskip\noindent
{\sectionsize\sectiontype \actualnumber.\ #1}
    \EQNcount=1
    \CLAIMcount=1
    \SUBSECTIONcount=1
    \nobreak\sectionskip\noindent}
\def\SECTIONNONR#1{\vskip0pt plus.3\vsize\penalty-75
    \vskip0pt plus -.3\vsize
    \global\advance\SECTIONcount by 1
    \beforesectionskip\noindent
{\sectionsize\sectiontype  #1}
     \EQNcount=1
     \CLAIMcount=1
     \SUBSECTIONcount=1
     \nobreak\sectionskip\noindent}
\def\SUBSECTION#1\par{\vskip0pt plus.2\vsize\penalty-75%
    \vskip0pt plus -.2\vsize%
    \beforesectionskip\noindent%
{\subsectionsize\subsectiontype \actualnumber.\number\SUBSECTIONcount.\ #1}
    \global\advance\SUBSECTIONcount by 1
    \nobreak\sectionskip\noindent}
\def\SUBSECTIONNONR#1\par{\vskip0pt plus.2\vsize\penalty-75
    \vskip0pt plus -.2\vsize
\beforesectionskip\noindent
{\subsectionsize\subsectiontype #1}
    \nobreak\sectionskip\noindent\noindent}
\def\ifundefined#1{\expandafter\ifx\csname#1\endcsname\relax}
\def\equ#1{\ifundefined{e#1}$\spadesuit$#1\else\csname e#1\endcsname\fi}
\def\clm#1{\ifundefined{c#1}$\spadesuit$#1\else\csname c#1\endcsname\fi}
\def\sec#1{\ifundefined{s#1}$\spadesuit$#1
\else Section \csname s#1\endcsname\fi}
\def\lab#1#2{\ifundefined{#1#2}$\spadesuit$#2\else\csname #1#2\endcsname\fi}
\def\fig#1{\ifundefined{fig#1}$\spadesuit$#1\else\csname fig#1\endcsname\fi}
\let\endarg=\par
\def\finish{\def\endarg{\par\endgroup}}
\def\start{\endarg\begingroup}

 \def\beginFROM{\start\parskip=0pt\vskip\baselineskip
\def\finish{\def\endarg{\egroup\par\endgroup}}
  \vbox\bgroup\obeylines\eightpoint\em\finish}

\def\ABSTRACT#1\par{
\vskip 1in {\noindent\sectionsize\sectiontype Abstract.} #1 \par}

\def\TODAY{\number\day~\ifcase\month\or January \or February \or March \or
April \or May \or June
\or July \or August \or September \or October \or November \or December \fi
\number\year\timecount=\number\time
\divide\timecount by 60
}
\newcount\timecount
\def\DRAFT{\def\lmargin##1{\strut\vadjust{\kern-\strutdepth
\vtop to \strutdepth{
\baselineskip\strutdepth\vss\rlap{\kern-1.2 truecm\eightpoint{##1}}}}}
\font\footfont=cmti7
\footline={{\footfont \hfil File:\jobname, \TODAY,  \number\timecount h}}
}
\newbox\strutboxJPE
\setbox\strutboxJPE=\hbox{\strut}
\def\subitem#1#2\par{\vskip\baselineskip\vskip-\ht\strutboxJPE{\item{#1}#2}}
\gdef\strutdepth{\dp\strutbox}
\def\lmargin#1{}
\def\hexnumber#1{\ifcase#1 0\or1\or2\or3\or4\or5\or6\or7\or8\or9\or
 A\or B\or C\or D\or E\or F\fi}
\textfont\msbfam=\tenmsb
\scriptfont\msbfam=\sevenmsb
\scriptscriptfont\msbfam=\fivemsb
\mathchardef\varkappa="0\hexnumber\msbfam7B%
\newcount\FIGUREcount \FIGUREcount=0
\newdimen\figcenter
\def\definefigure#1{\global\advance\FIGUREcount by 1%
\NEWDEF{fig}{#1}{Fig.\ \number\FIGUREcount}
\immediate\write16{  FIG \number\FIGUREcount : #1}}
\def\figure#1#2#3#4\cr{\null%
\definefigure{#1}
{\goodbreak\figcenter=\hsize\relax
\advance\figcenter by -#3truecm
\divide\figcenter by 2
\midinsert\vskip #2truecm\noindent\hskip\figcenter
\includegraphics{#1}\vskip 0.8truecm\noindent \vbox{\eightpoint\noindent
{\bf\fig{#1}}: #4}\endinsert}}
\def\figurewithtex#1#2#3#4#5\cr{\null%
\definefigure{#1}
{\goodbreak\figcenter=\hsize\relax
\advance\figcenter by -#4truecm
\divide\figcenter by 2
\midinsert\vskip #3truecm\noindent\hskip\figcenter
\includegraphics{#1}{\hskip\texpscorrection\input #2 }\vskip 0.8truecm\noindent \vbox{\eightpoint\noindent
{\bf\fig{#1}}: #5}\endinsert}}
\def\figurewithtexplus#1#2#3#4#5#6\cr{\null%
\definefigure{#1}
{\goodbreak\figcenter=\hsize\relax
\advance\figcenter by -#4truecm
\divide\figcenter by 2
\midinsert\vskip #3truecm\noindent\hskip\figcenter
\includegraphics{#1}{\hskip\texpscorrection\input #2 }\vskip #5truecm\noindent \vbox{\eightpoint\noindent
{\bf\fig{#1}}: #6}\endinsert}}
\catcode`@=11
\def\footnote#1{\let\@sf\empty 
  \ifhmode\edef\@sf{\spacefactor\the\spacefactor}\/\fi
  #1\@sf\vfootnote{#1}}
\def\vfootnote#1{\insert\footins\bgroup\eightpoint
  \interlinepenalty\interfootnotelinepenalty
  \splittopskip\ht\strutbox 
  \splitmaxdepth\dp\strutbox \floatingpenalty\@MM
  \leftskip\z@skip \rightskip\z@skip \spaceskip\z@skip \xspaceskip\z@skip
  \textindent{#1}\footstrut\futurelet\next\fo@t}
\def\fo@t{\ifcat\bgroup\noexpand\next \let\next\f@@t
  \else\let\next\f@t\fi \next}
\def\f@@t{\bgroup\aftergroup\@foot\let\next}
\def\f@t#1{#1\@foot}
\def\@foot{\strut\egroup}
\def\footstrut{\vbox to\splittopskip{}}
\skip\footins=\bigskipamount 
\count\footins=1000 
\dimen\footins=8in 
\catcode`@=12 

\def\AA{{\script A}}
\def\BB{{\script B}}
\def\CC{{\script C}}

\def\HH{{\script H}}
\def\LL{{\script L}}

\def\OO{{\script O}}
\def\PP{{\script P}}

\def\HALF{{\textstyle{1\over 2}}}

\def\QED{\hfill\smallskip
         \line{$\hfill{\vcenter{\vbox{\hrule height 0.2pt
	\hbox{\vrule width 0.2pt height 1.8ex \kern 1.8ex
		\vrule width 0.2pt}
	\hrule height 0.2pt}}}$
               \ \ \ \ \ \ }
         \bigskip}
\def\real{{\bf R}}

\def\integer{{\bf Z}}

\def\PROOF{\medskip\noindent{\bf Proof.\ }}
\def\REMARK{\medskip\noindent{\bf Remark.\ }}
\def\LIKEREMARK#1{\medskip\noindent{\bf #1.\ }}
\normalbaselineskip=5.25mm
\baselineskip=5.25mm
\parskip=10pt
\beforesectionskipamount=24pt plus8pt minus8pt
\sectionskipamount=3pt plus1pt minus1pt
\def\em{\it}
\tenpoint
\null
\catcode`\@=11
\ifx\amstexloaded@\relax\catcode`\@=\active
\endinput\fi
\catcode`\@=\active
\def\period{\unskip.\spacefactor3000 { }}
%
%
\newbox\noboxJPE
\newbox\byboxJPE
\newbox\paperboxJPE
\newbox\yrboxJPE
\newbox\jourboxJPE
\newbox\pagesboxJPE
\newbox\volboxJPE
\newbox\preprintboxJPE
\newbox\toappearboxJPE
\newbox\bookboxJPE
\newbox\bybookboxJPE
\newbox\publisherboxJPE
\newbox\inprintboxJPE
\def\refclearJPE{
   \setbox\noboxJPE=\null             \gdef\isnoJPE{F}
   \setbox\byboxJPE=\null             \gdef\isbyJPE{F}
   \setbox\paperboxJPE=\null          \gdef\ispaperJPE{F}
   \setbox\yrboxJPE=\null             \gdef\isyrJPE{F}
   \setbox\jourboxJPE=\null           \gdef\isjourJPE{F}
   \setbox\pagesboxJPE=\null          \gdef\ispagesJPE{F}
   \setbox\volboxJPE=\null            \gdef\isvolJPE{F}
   \setbox\preprintboxJPE=\null       \gdef\ispreprintJPE{F}
   \setbox\toappearboxJPE=\null       \gdef\istoappearJPE{F}
   \setbox\inprintboxJPE=\null        \gdef\isinprintJPE{F}
   \setbox\bookboxJPE=\null           \gdef\isbookJPE{F}  \gdef\isinbookJPE{F}
     
   \setbox\bybookboxJPE=\null         \gdef\isbybookJPE{F}
   \setbox\publisherboxJPE=\null      \gdef\ispublisherJPE{F}
}

\def\ref{\refclearJPE}
\def\no#1{\gdef\isnoJPE{T}\setbox\noboxJPE=\hbox{#1}}
\def\by#1{\gdef\isbyJPE{T}\setbox\byboxJPE=\hbox{#1}}
\def\paper#1{\gdef\ispaperJPE{T}\setbox\paperboxJPE=\hbox{#1}}
\def\yr#1{\gdef\isyrJPE{T}\setbox\yrboxJPE=\hbox{#1}}
\def\jour#1{\gdef\isjourJPE{T}\setbox\jourboxJPE=\hbox{#1}}
\def\pages#1{\gdef\ispagesJPE{T}\setbox\pagesboxJPE=\hbox{#1}}
\def\vol#1{\gdef\isvolJPE{T}\setbox\volboxJPE=\hbox{\bf #1}}
\def\preprint#1{\gdef
\ispreprintJPE{T}\setbox\preprintboxJPE=\hbox{#1}}

\def\book#1{\gdef\isbookJPE{T}\setbox\bookboxJPE=\hbox{\em #1}}
\def\publisher#1{\gdef
\ispublisherJPE{T}\setbox\publisherboxJPE=\hbox{#1}}
\def\inbook#1{\gdef\isinbookJPE{T}\setbox\bookboxJPE=\hbox{\em #1}}
\def\bybook#1{\gdef\isbybookJPE{T}\setbox\bybookboxJPE=\hbox{#1}}
\newdimen\refindent
\refindent=5em
\def\endref{\sfcode`.=1000
 \if T\isnoJPE
 \hangindent\refindent\hangafter=1
      \noindent\hbox to\refindent{[\unhbox\noboxJPE\unskip]\hss}\ignorespaces
     \else  \noindent    \fi
 \if T\isbyJPE    \unhbox\byboxJPE\unskip: \fi
 \if T\ispaperJPE \unhbox\paperboxJPE\unskip\period \fi
 \if T\isbookJPE {\it\unhbox\bookboxJPE\unskip}\if T\ispublisherJPE, \else.
\fi\fi
 \if T\isinbookJPE In {\it\unhbox\bookboxJPE\unskip}\if T\isbybookJPE,
\else\period \fi\fi
 \if T\isbybookJPE  (\unhbox\bybookboxJPE\unskip)\period \fi
 \if T\ispublisherJPE \unhbox\publisherboxJPE\unskip \if T\isjourJPE, \else\if
T\isyrJPE \  \else\period \fi\fi\fi
 \if T\istoappearJPE (To appear)\period \fi
 \if T\ispreprintJPE Pre\-print\period \fi
 \if T\isjourJPE    \unhbox\jourboxJPE\unskip\ \fi
 \if T\isvolJPE     \unhbox\volboxJPE\unskip\if T\ispagesJPE, \else\ \fi\fi
 \if T\ispagesJPE   \unhbox\pagesboxJPE\unskip\  \fi
 \if T\isyrJPE      (\unhbox\yrboxJPE\unskip)\period \fi
 \if T\isinprintJPE (in print)\period \fi
\filbreak
}
\normalbaselineskip=12pt
\baselineskip=12pt
\parskip=0pt
\parindent=22.222pt
\beforesectionskipamount=24pt plus0pt minus6pt
\sectionskipamount=7pt plus3pt minus0pt
\overfullrule=0pt
\hfuzz=2pt
\nopagenumbers
\headline={\ifnum\pageno>1 {\hss\tenrm-\ \folio\ -\hss} \else
{\hfill}\fi}
\if F\PSfonts

\font\toplinefont=cmr10
\font\pagenumberfont=cmr10
\let\tenpoint=\rm
\else

\font\toplinefont=cmcsc10
\font\pagenumberfont=ptmb at 10pt
\fi
\newdimen\itemindent\itemindent=1.5em

\def\textindent#1{\indent\llap{#1\enspace}\ignorespaces}
\def\item{\par\noindent
\hangindent\itemindent\hangafter=1\relax
\setitemmark}
\def\setitemindent#1{\setbox0=\hbox{\ignorespaces#1\unskip\enspace}%
\itemindent=\wd0\relax
\message{|\string\setitemindent: Mark width modified to hold
         |`\string#1' plus an \string\enspace\space gap. }%
}
\def\setitemmark#1{\checkitemmark{#1}%
\hbox to\itemindent{\hss#1\enspace}\ignorespaces}
\def\checkitemmark#1{\setbox0=\hbox{\enspace#1}%
\ifdim\wd0>\itemindent
   \message{|\string\item: Your mark `\string#1' is too wide. }%
\fi}
\def\SECTION#1{\vskip0pt plus.2\vsize\penalty-75
    \vskip0pt plus -.2\vsize
    \global\advance\SECTIONcount by 1
    \beforesectionskip\noindent
{\sectionsize\sectiontype \actualnumber.\ #1}
    \EQNcount=1
    \CLAIMcount=1
    \SUBSECTIONcount=1
    \nobreak\sectionskip\noindent}
\font\fteenrm=ptmr at 14pt
\font\fteeni=cmmi10 at 14pt
\font\fteenbf=ptmb at 14pt
\font\fteensy=cmsy10 at	14pt
\font\fteenex=cmex10	at	14pt
\font\fteenit=cmti10 at	14pt
\font\fteenmsb=msbm10	at	14pt
\font\fteentruecmsy=cmsy10	at	14pt
\font\fteentruecmr=cmr10	at	14pt

\def\fteenpoint{\def\rm{\fam0\fteenrm}%
\textfont0=\fteenrm
  \scriptfont0=\twelverm
  \scriptscriptfont0=\tenrm 
\textfont1=\fteeni
  \scriptfont1=\twelvei
  \scriptscriptfont1=\teni 
\textfont2=\fteensy
  \scriptfont2=\twelvesy
  \scriptscriptfont2=\tensy 
\textfont3=\fteenex
  \scriptfont3=\fteenex
  \scriptscriptfont3=\fteenex 
\textfont\itfam=\fteenit
  \scriptfont\itfam=\twelveit
  \scriptscriptfont\itfam=\tenit 
  \def\it{\fam\itfam\fteenit}%
\textfont\bffam=\fteenbf
  \scriptfont\bffam=\twelvebf
  \scriptscriptfont\bffam=\tenbf
  \def\bf{\fam\bffam\fteenbf}%
\textfont\msbfam=\fteenmsb
  \scriptfont\msbfam=\twelvemsb
  \scriptscriptfont\msbfam=\tenmsb
\textfont\truecmr=\fteentruecmr
  \scriptfont\truecmr=\twelvetruecmr
  \scriptscriptfont\truecmr=\tentruecmr
  \def\truerm{\fam\truecmr\fteentruecmr}%
\textfont\truecmsy=\fteentruecmsy
  \scriptfont\truecmsy=\twelvetruecmsy
  \scriptscriptfont\truecmsy=\tentruecmsy
\setbox\strutbox=\hbox{\vrule height7pt depth2pt width0pt}%
\normalbaselineskip=16.5pt
\normalbaselines
\rm
}

\expandafter\xdef\csname
sintroduction\endcsname{1}
\expandafter\xdef\csname
estart1\endcsname{(1.1)}
\expandafter\xdef\csname
e11a\endcsname{(1.2)}
\expandafter\xdef\csname
eparabolic\endcsname{(1.3)}
\expandafter\xdef\csname
ehyperbolic\endcsname{(1.4)}
\expandafter\xdef\csname
ebbbb\endcsname{(1.5)}
\expandafter\xdef\csname
scoercive\endcsname{2}
\expandafter\xdef\csname
eproblem\endcsname{(2.1)}
\expandafter\xdef\csname
ehh\endcsname{(2.2)}
\expandafter\xdef\csname
egdd\endcsname{(2.3)}
\expandafter\xdef\csname
ealphaloc1\endcsname{(2.4)}
\expandafter\xdef\csname
ebb1def\endcsname{(2.5)}
\expandafter\xdef\csname
cbbloc1\endcsname{Definition\ 2.1}
\expandafter\xdef\csname
eF00\endcsname{(2.6)}
\expandafter\xdef\csname
eG0\endcsname{(2.7)}
\expandafter\xdef\csname
egbound\endcsname{(2.8)}
\expandafter\xdef\csname
edtE00\endcsname{(2.9)}
\expandafter\xdef\csname
cE0\endcsname{Lemma\ 2.2}
\expandafter\xdef\csname
wg\endcsname{C_{0}}
\expandafter\xdef\csname
wg1\endcsname{C_{1}}
\expandafter\xdef\csname
wg2\endcsname{C_{2}}
\expandafter\xdef\csname
edtE0\endcsname{(2.10)}
\expandafter\xdef\csname
cE00\endcsname{Lemma\ 2.3}
\expandafter\xdef\csname
wh\endcsname{C_{3}}
\expandafter\xdef\csname
ethisisit\endcsname{(2.11)}
\expandafter\xdef\csname
clinfty\endcsname{Proposition\ 2.4}
\expandafter\xdef\csname
eaux\endcsname{(2.12)}
\expandafter\xdef\csname
wSobolev\endcsname{C_{4}}
\expandafter\xdef\csname
eSobolev\endcsname{(2.13)}
\expandafter\xdef\csname
cSobolev\endcsname{Lemma\ 2.5}
\expandafter\xdef\csname
eubound\endcsname{(2.14)}
\expandafter\xdef\csname
ewz\endcsname{(2.15)}
\expandafter\xdef\csname
eF10\endcsname{(2.16)}
\expandafter\xdef\csname
elong\endcsname{(2.17)}
\expandafter\xdef\csname
wf2\endcsname{C_{5}}
\expandafter\xdef\csname
elong2\endcsname{(2.18)}
\expandafter\xdef\csname
w1\endcsname{C_{6}}
\expandafter\xdef\csname
elong3\endcsname{(2.19)}
\expandafter\xdef\csname
w777\endcsname{C_{7}}
\expandafter\xdef\csname
w7777\endcsname{C_{8}}
\expandafter\xdef\csname
w77\endcsname{C_{9}}
\expandafter\xdef\csname
w7\endcsname{C_{10}}
\expandafter\xdef\csname
egood0\endcsname{(2.20)}
\expandafter\xdef\csname
egood\endcsname{(2.21)}
\expandafter\xdef\csname
cf1\endcsname{Theorem\ 2.6}
\expandafter\xdef\csname
egood3\endcsname{(2.22)}
\expandafter\xdef\csname
egood4\endcsname{(2.23)}
\expandafter\xdef\csname
slinearized\endcsname{3}
\expandafter\xdef\csname
ediff1\endcsname{(3.1)}
\expandafter\xdef\csname
ellform\endcsname{(3.2)}
\expandafter\xdef\csname
ehfirst\endcsname{(3.3)}
\expandafter\xdef\csname
everylong\endcsname{(3.4)}
\expandafter\xdef\csname
wH\endcsname{C_{11}}
\expandafter\xdef\csname
eHfinal\endcsname{(3.5)}
\expandafter\xdef\csname
eldd\endcsname{(3.6)}
\expandafter\xdef\csname
eh2def\endcsname{(3.7)}
\expandafter\xdef\csname
cH2\endcsname{Definition\ 3.1}
\expandafter\xdef\csname
ebb2def\endcsname{(3.8)}
\expandafter\xdef\csname
cbbloc2\endcsname{Definition\ 3.2}
\expandafter\xdef\csname
egrowth\endcsname{(3.9)}
\expandafter\xdef\csname
cfrakuv\endcsname{Theorem\ 3.3}
\expandafter\xdef\csname
sprep\endcsname{4}
\expandafter\xdef\csname
eQdef\endcsname{(4.1)}
\expandafter\xdef\csname
ehdef\endcsname{(4.2)}
\expandafter\xdef\csname
wM\endcsname{C_{12}}
\expandafter\xdef\csname
eboundona\endcsname{(4.3)}
\expandafter\xdef\csname
cq\endcsname{Lemma\ 4.1}
\expandafter\xdef\csname
ekernel\endcsname{(4.4)}
\expandafter\xdef\csname
wQ<\endcsname{C_{13}}
\expandafter\xdef\csname
eQdef0\endcsname{(4.5)}
\expandafter\xdef\csname
eboundonq\endcsname{(4.6)}
\expandafter\xdef\csname
cQa\endcsname{Corollary\ 4.2}
\expandafter\xdef\csname
shigh\endcsname{5}
\expandafter\xdef\csname
ekaa\endcsname{(5.1)}
\expandafter\xdef\csname
wnu0\endcsname{C_{14}}
\expandafter\xdef\csname
enumore\endcsname{(5.2)}
\expandafter\xdef\csname
ckaa\endcsname{Proposition\ 5.1}
\expandafter\xdef\csname
wnu\endcsname{C_{15}}
\expandafter\xdef\csname
enufinal\endcsname{(5.3)}
\expandafter\xdef\csname
shighfreq\endcsname{6}
\expandafter\xdef\csname
eQdef2\endcsname{(6.1)}
\expandafter\xdef\csname
epdef\endcsname{(6.2)}
\expandafter\xdef\csname
wQ\endcsname{C_{16}}
\expandafter\xdef\csname
ePbound\endcsname{(6.3)}
\expandafter\xdef\csname
cq2\endcsname{Lemma\ 6.1}
\expandafter\xdef\csname
wQQ\endcsname{C_{17}}
\expandafter\xdef\csname
cQonB\endcsname{Lemma\ 6.2}
\expandafter\xdef\csname
ejdef0\endcsname{(6.4)}
\expandafter\xdef\csname
ehdd\endcsname{(6.5)}
\expandafter\xdef\csname
wdd\endcsname{C_{18}}
\expandafter\xdef\csname
ehdd2\endcsname{(6.6)}
\expandafter\xdef\csname
esplit\endcsname{(6.7)}
\expandafter\xdef\csname
eeta0\endcsname{(6.8)}
\expandafter\xdef\csname
edefk*\endcsname{(6.9)}
\expandafter\xdef\csname
eddx0\endcsname{(6.10)}
\expandafter\xdef\csname
edeltadef\endcsname{(6.11)}
\expandafter\xdef\csname
euvbound\endcsname{(6.12)}
\expandafter\xdef\csname
ehdd3\endcsname{(6.13)}
\expandafter\xdef\csname
ehdd3a\endcsname{(6.14)}
\expandafter\xdef\csname
enubounds\endcsname{(6.15)}
\expandafter\xdef\csname
ehdd4\endcsname{(6.16)}
\expandafter\xdef\csname
eenfin\endcsname{(6.17)}
\expandafter\xdef\csname
ejdef2\endcsname{(6.18)}
\expandafter\xdef\csname
egammadef\endcsname{(6.19)}
\expandafter\xdef\csname
cJ\endcsname{Proposition\ 6.3}
\expandafter\xdef\csname
cuv\endcsname{Lemma\ 6.4}
\expandafter\xdef\csname
eJL\endcsname{(6.20)}
\expandafter\xdef\csname
eLt\endcsname{(6.21)}
\expandafter\xdef\csname
cJL\endcsname{Theorem\ 6.5}
\expandafter\xdef\csname
sfull\endcsname{7}
\expandafter\xdef\csname
ediff1a\endcsname{(7.1)}
\expandafter\xdef\csname
ellforma\endcsname{(7.2)}
\expandafter\xdef\csname
emmm\endcsname{(7.3)}
\expandafter\xdef\csname
erepr1\endcsname{(7.4)}
\expandafter\xdef\csname
wMM1\endcsname{C_{19}}
\expandafter\xdef\csname
wMM2\endcsname{C_{20}}
\expandafter\xdef\csname
eMM1\endcsname{(7.5)}
\expandafter\xdef\csname
cexpo1\endcsname{Proposition\ 7.1}
\expandafter\xdef\csname
ellformb\endcsname{(7.6)}
\expandafter\xdef\csname
elook\endcsname{(7.7)}
\expandafter\xdef\csname
wMM9\endcsname{C_{21}}
\expandafter\xdef\csname
wMM3\endcsname{C_{22}}
\expandafter\xdef\csname
wMM4\endcsname{C_{23}}
\expandafter\xdef\csname
eMM2\endcsname{(7.8)}
\expandafter\xdef\csname
cexpo2\endcsname{Proposition\ 7.2}
\expandafter\xdef\csname
sbounds\endcsname{8}
\expandafter\xdef\csname
enewnorm\endcsname{(8.1)}
\expandafter\xdef\csname
wball\endcsname{C_{24}}
\expandafter\xdef\csname
eballs\endcsname{(8.2)}
\expandafter\xdef\csname
cballs\endcsname{Theorem\ 8.1}
\expandafter\xdef\csname
etaudef\endcsname{(8.3)}
\expandafter\xdef\csname
etheineq\endcsname{(8.4)}
\expandafter\xdef\csname
ebound1\endcsname{(8.5)}
\expandafter\xdef\csname
winitial\endcsname{C_{25}}
\expandafter\xdef\csname
ebound2\endcsname{(8.6)}
\expandafter\xdef\csname
edecomp1\endcsname{(8.7)}
\expandafter\xdef\csname
efinaleta0\endcsname{(8.8)}
\expandafter\xdef\csname
efinal1\endcsname{(8.9)}
\expandafter\xdef\csname
ecart\endcsname{(8.10)}
\expandafter\xdef\csname
wBB\endcsname{C_{26}}
\expandafter\xdef\csname
eincomplete\endcsname{(8.11)}
\expandafter\xdef\csname
czzz\endcsname{Lemma\ 8.2}
\expandafter\xdef\csname
eexponential\endcsname{(8.12)}
\expandafter\xdef\csname
wbound\endcsname{C_{27}}
\expandafter\xdef\csname
ethebound\endcsname{(8.13)}
\expandafter\xdef\csname
ccart\endcsname{Lemma\ 8.3}
\expandafter\xdef\csname
wx4\endcsname{C_{28}}
\expandafter\xdef\csname
wx5\endcsname{C_{29}}
\expandafter\xdef\csname
eder\endcsname{(8.14)}
\expandafter\xdef\csname
wx44\endcsname{C_{30}}
\expandafter\xdef\csname
wx99\endcsname{C_{31}}
\expandafter\xdef\csname
wh22\endcsname{C_{32}}
\expandafter\xdef\csname
eS\endcsname{(8.15)}
\expandafter\xdef\csname
esh1\endcsname{(8.16)}
\expandafter\xdef\csname
wn1\endcsname{C_{33}}
\expandafter\xdef\csname
wL\endcsname{C_{34}}
\expandafter\xdef\csname
efor9\endcsname{(8.17)}
\expandafter\xdef\csname
epart1\endcsname{(8.18)}
\expandafter\xdef\csname
wL2\endcsname{C_{35}}
\expandafter\xdef\csname
cfinite\endcsname{Proposition\ 8.4}
\expandafter\xdef\csname
wp1\endcsname{C_{36}}
\expandafter\xdef\csname
wp2\endcsname{C_{37}}
\expandafter\xdef\csname
wdelta\endcsname{C_{38}}
\expandafter\xdef\csname
wdelta2\endcsname{C_{39}}
\expandafter\xdef\csname
sentropy\endcsname{9}
\expandafter\xdef\csname
ennn\endcsname{(9.1)}
\expandafter\xdef\csname
ewwwdef\endcsname{(9.2)}
\expandafter\xdef\csname
ssubadditivity\endcsname{9.2}
\expandafter\xdef\csname
esmooth\endcsname{(9.3)}
\expandafter\xdef\csname
cconnect\endcsname{Theorem\ 9.1}
\expandafter\xdef\csname
cboth\endcsname{Corollary\ 9.2}
\expandafter\xdef\csname
elemme1\endcsname{(9.4)}
\expandafter\xdef\csname
eboundary\endcsname{(9.5)}
\expandafter\xdef\csname
egdef\endcsname{(9.6)}
\expandafter\xdef\csname
erand\endcsname{(9.7)}
\expandafter\xdef\csname
capprox\endcsname{Proposition\ 9.3}
\expandafter\xdef\csname
er\endcsname{(9.8)}
\expandafter\xdef\csname
erightend\endcsname{(9.9)}
\expandafter\xdef\csname
enumu\endcsname{(9.10)}
\expandafter\xdef\csname
conestep\endcsname{Lemma\ 9.4}
\expandafter\xdef\csname
ejdef\endcsname{(9.11)}
\expandafter\xdef\csname
efg\endcsname{(9.12)}
\expandafter\xdef\csname
efg2\endcsname{(9.13)}
\expandafter\xdef\csname
eupper\endcsname{(9.14)}
\expandafter\xdef\csname
elower\endcsname{(9.15)}
\expandafter\xdef\csname
efg3\endcsname{(9.16)}
\expandafter\xdef\csname
egood1\endcsname{(9.17)}
\expandafter\xdef\csname
erho\endcsname{(9.18)}
\expandafter\xdef\csname
egood2\endcsname{(9.19)}
\expandafter\xdef\csname
see\endcsname{9.3}
\expandafter\xdef\csname
eballs2\endcsname{(9.20)}
\expandafter\xdef\csname
eballs5\endcsname{(9.21)}
\expandafter\xdef\csname
ecompare\endcsname{(9.22)}
\expandafter\xdef\csname
esub\endcsname{(9.23)}
\expandafter\xdef\csname
w6\endcsname{C_{40}}
\expandafter\xdef\csname
w66\endcsname{C_{41}}
\expandafter\xdef\csname
wee\endcsname{C_{42}}
\expandafter\xdef\csname
eee\endcsname{(9.24)}
\expandafter\xdef\csname
cee\endcsname{Theorem\ 9.5}
\expandafter\xdef\csname
stop\endcsname{9.4}
\expandafter\xdef\csname
ctoto\endcsname{Lemma\ 9.6}
\expandafter\xdef\csname
elim\endcsname{(9.25)}
\expandafter\xdef\csname
cexist\endcsname{Theorem\ 9.7}
\expandafter\xdef\csname
eregion\endcsname{(9.26)}
\expandafter\xdef\csname
wlast\endcsname{C_{43}}

\setitemindent{iii)}
\def\frac#1#2{{#1\over #2}}
\def\KK{{\cal K}}
\def\attra{\GG}
\def\calu{{\cal U}}
\def\idx{\int}

\def\supp{{\rm{supp\ }}}

\def\dx{{\kern -0.2em{\rm d}x\,}}
\def\d#1{{\kern -0.2em{\rm d}#1\,}}
\def\dpp{{\kern -0.2em{\rm d}p\,}}
\def\GG{{\script G}}
\def\QQ{{\script Q}}
\def\SS{{\script S}}
\def\L{{\rm L}}

\def\fff{\varphi}
\def\ggg{\psi}
\def\hdd{h_\delta}
\def\hh{h}
\def\hhh{h}

\def\loc{{\rm loc}}
\def\V#1#2{\left({{#1}\atop {#2}}\right)}
\def\FOUR{{\textstyle{1\over 4}}}
\def\EIGHT{{\textstyle{1\over 8}}}
\def\epsilonp{\alpha}
\let\epsilon=\varepsilon
\let\theta=\vartheta
\let\kappa=\varkappa
\def\kstar{k_*}
\let\phi=\varphi
\def\fu{{\frak u}}
\def\fv{{\frak v}}

\def\gu{\fu}
\def\gv{\fv}
\def\nut#1{{}}
\def\FABC{F_{\epsilon ,A,B,C}}
\def\FABCnull{F^0_{\delta ,A,B,C,\xi}}
\def\E{E_{\epsilon,G,g_\l,g_\r}}
\def\WWW{W^{1,\infty }}
\def\FF{{\cal F}}
\def\CC{{\cal C}}
\def\GG{{\cal G}}
\def\R{{R }}
\def\l{{\rm L}}
\def\r{{\rm R}}
\def\cartarg{{\sigma x\over 2}- {\pi j\over 4}}
\def\cartargfour{{2\sigma x}- {\pi j}}
\def\cartk{k_* x - {\pi j\over 4}}
\def\lpoint{2L\kstar}
\newcount\Xvalue\Xvalue=10
\newcount\Tvalue\Tvalue=10
\def\X{\the\Xvalue}
\def\T{\the\Tvalue}
\newcount\TXvalue\TXvalue=\the\Xvalue
\multiply\TXvalue by \the \Tvalue
\newcount\TNUvalue\TNUvalue=40
\newcount\TMUvalue\TMUvalue=\the\TXvalue
\multiply\TMUvalue by 2
\newcount\TRvalue\TRvalue=40
\newcount\TRplusvalue\TRplusvalue=\the\TRvalue
\advance\TRplusvalue by 1
\newcount\Ttemp
\def\TX{\the\TXvalue}
\def\TNU{\the\TNUvalue}
\def\TMU{\the\TMUvalue}
\def\TR{\the\TRvalue}
\def\wqeps{{\cal W}_{I}^{\epsilon}}
\def\text#1{\leavevmode\hbox{#1}}
\let\truett=\tt
\fontdimen3\tentt=2pt\fontdimen4\tentt=2pt
\def\tt{\hfill\break\null\kern -2truecm\truett ************ }
\newcount\Ccount\Ccount=-1
\def\C#1{\lmargin{#1}\global\advance\Ccount by 1
\NEWDEF{w}{#1}{C_{\number\Ccount}}
C_{\number\Ccount}
}
\def\CX#1{\ifundefined{w#1}\spadesuit#1\else\csname w#1\endcsname\fi}
{\fteenpoint{\centerline{\bf Topological Entropy and $\varepsilon $-Entropy }}}
\vskip 0.5truecm
{\fteenpoint{\centerline{{\bf for Damped Hyperbolic Equations}}}
\vskip 0.5truecm}
{\it{\centerline{ P. Collet${}^{1}$ and J.-P. Eckmann${}^{2,3}$}}}
\vskip 0.3truecm
{\eightpoint
\centerline{${}^1$Centre de Physique Th\'eorique, Laboratoire CNRS UMR
7644,}\centerline{
Ecole Polytechnique, F-91128 Palaiseau Cedex, France}
\centerline{${}^2$D\'ept.~de Physique Th\'eorique, Universit\'e de Gen\`eve,
CH-1211 Gen\`eve 4, Switzerland}
\centerline{${}^3$Section de Math\'ematiques, Universit\'e de Gen\`eve,
CH-1211 Gen\`eve 4, Switzerland}
}
\vskip 0.5truecm\headline
{\ifnum\pageno>1 {\toplinefont Damped Hyperbolic Equations}
\hfill{\pagenumberfont\folio}\fi}
\vskip 1cm

\narrower
{\eightpoint\baselineskip 11pt
\LIKEREMARK{Abstract}{We study damped hyperbolic equations on the
infinite line. We show that on the global attracting set $\GG$ the
$\epsilon$-entropy (per unit length) exists in the topology of
$\WWW$. We also show that the topological entropy per unit length
of $\GG$ exists.
These results are shown using two main techniques: Bounds in bounded
domains
in position space and for large momenta, and a novel submultiplicativity
argument in $\WWW$.
}}

\advance\leftskip-\parindent
\advance\rightskip-\parindent
\SECT{introduction}{Introduction}

This paper is an extension of our earlier papers [CE1, CE2] to mixed
parabolic-hyperbolic equations in the infinite domain $\real$:
$$
\eta^2\partial_t^2 u(x,t)+\partial _t u(x,t)\,=\, 
\partial_x^2 u(x,t) +U'\bigl(u(x,t)\bigr)~,
\EQ{start1} 
$$
where $U(s)=s^2/2-s^4/4$. The particular choice of the potential $U$
is in fact not very important, but we will deal only with this one.
This problem can be written as a system:
$$
\eqalign{
\partial_t u(x,t)\,&=\,v(x,t)~,\cr
\eta^2\partial_t v(x,t) \,&=\,-v(x,t)+\partial_x^2 u(x,t) +U'\bigl(u(x,t)\bigr)~.\cr
}\EQ{11a}
$$
The functions $u$ will be real-valued, but the extension to
vector-valued functions is easy and is left to the reader, since it
only complicates notation.

The question we ask is about the nature of the attracting set for this
problem, its complexity, and in particular its $\epsilon $-entropy. We
have developed this subject for parabolic problems in the two papers
described above and we study now 
the complexity in this parabolic-hyperbolic setting. The difference
with the parabolic case is the absence of regularization. In the
parabolic case, the dispersion law is, written in Fourier space for
the linearized equation
$$
\partial_t \tilde u(k,t)\,=\,(1-k^2 )\tilde u(k,t)~,
\EQ{parabolic} 
$$
when $U(s)=s^2/2+\OO(s^3)$ near $s=0$. In the case we consider now,
the problem is rather a system of the form
$$
\eqalign{
\partial _t \tilde u(k,t)\,&=\,\tilde v(k,t)~,\cr
\eta^2\partial _t\tilde v(k,t)\,&=\,-\tilde v(k,t)+(1-k^2)\tilde u(k,t)~.\cr
}
\EQ{hyperbolic} 
$$
Thus, as is well known, \equ{parabolic} regularizes derivatives
because
$|k|\exp\bigl((1-k^2)t\bigr)$ is bounded in $k$ when $t>0$, while the
real part of the eigenvalue of the system \equ{hyperbolic} is, for
large $|k|$, only as negative as $-\OO(\eta^{-2})$, and therefore
the exponential is only bounded like $|k|\exp\bigl(-C\eta^{-2}t\bigr)$
for some $C>0$. This diverges with $|k|$, but converges
(non-uniformly) to 0 as $t\to\infty $.

One can ask whether this reduced form of regularization manifests
itself in an increased complexity of either the attracting set, or
some forward invariant set of bounded initial data. The conclusion of
our paper is that {\em the complexity of the problem \equ{hyperbolic}
is of the same order as that of \equ{parabolic}.}

%
Since we work on the infinite line, we need local topologies. This
will be achieved by choosing a cutoff function $h$:
$$
h_\delta(x)\,=\,{1\over (1+\delta^2x^2)^2}~.
$$
We could take other functions with sufficiently strong
polynomial decay, but the nice ideas of Mielke [M1, M2] using
exponentially decaying cutoff functions do not seem to work here.
We then consider local Sobolev norms of the form
$$
\|(u,v)\|_{\hdd,2}^2\,=\,\int\dx \hdd(x) \bigg(u^2
+2(u')^2+(u'')^2+\eta^2\bigl(v^2+(v')^2\bigr)\bigg)(x) ~, 
$$
and then local spaces $\HH_{\delta,\loc}^2$ with the norm
$$
\|(u,v)\|_{\delta,\loc,2}\,=\,\sup_{\xi\in\real}\|(u,v)\|_{h_{\delta,\xi},2}~,
$$
where ${h_{\delta,\xi}}(x)=\hdd(x-\xi)$.
Note that this norm, and many others used in this paper, has {\em one
more derivative in $u$ than in $v$}. Such norms are typical when one
writes equations such as \equ{start1} with two components as in
\equ{11a}. 

We will show that every initial condition with finite
$\|(u,v)\|_{\delta,\loc,2}$ will end up after some finite time in a {\em
bounded} set in this norm. We call this bounded set the attracting
set $\GG$. The attractor $\AA$ is then defined as 
$$
\AA\,=\,\bigcap_{t>0} \Phi^t(\GG)~,
$$
where $t\mapsto\Phi^t$ is the flow defined by \equ{start1}. We will not
only study the complexity of $\AA$, but we can also make statements
about functions which have ``evolved for long enough,'' namely functions
in
$\GG_T\equiv\bigcap_{T>t>0} \Phi^t(\GG)$ for some large $T$.
Indeed, given some interval $[-L,L]$ in $\real$, with $L\gg 1/\epsilon
$, we show in \sec{ee}
that
$\GG_T$ {\em when restricted to $[-L,L]$ in the 
variable $x$} can be covered by $N_L(\epsilon )=C^{L\log{1/\epsilon
}}$ balls of 
radius $\epsilon $ in $\HH_{\delta,\loc}^2$.
Our argument does not rely on compactness, but only on a comparison of
the
number of balls with radius $\epsilon $ relative to the
number needed when the radius is $2\epsilon $:
$$
N_L(\epsilon )\,\le\,N_{L+A/\epsilon }(2\epsilon )C^L~,
$$
for some constants $C$ and $A$ (\clm{finite} and \equ{balls2}).
It is at this point that we use the invariance of $\GG$, the fact that
high-momentum parts of the solution are damped with an exponential
rate of about $\eta^{-2}$, and that the low momentum parts are Fourier
transforms of analytic functions, which can be finitely sampled by the
Cartwright formula \equ{cart}.

We then change topology to $\WWW$ (functions in $\L^\infty $ whose
derivatives are also in $\L^\infty $) and show that the results
obtained for $\HH_{\delta,\loc}^2$ give bounds in $\WWW$. We introduce a
new type of submultiplicativity bound in \sec{subadditivity}. Indeed, we
show in \clm{both} (which is an easy consequence of the \clm{connect})
that if a bounded set of functions in $\CC^2$ 
can be covered by $N_1$ balls of radius
$\epsilon $ in $\WWW_{I_1}$ and by $N_2$ on $\WWW_{I_2}$, where $I_1$
and $I_2$ are disjoint intervals, then it can be covered by
$$
C(\epsilon )N_1 N_2 
\EQ{bbbb}
$$
balls of radius $\epsilon $ in $\WWW_{I_1\cup I_2}$. The point here is
that $C$ only depends on $\epsilon $ (and the bound on the functions)
and that the balls on $I_1\cup
I_2$ have the {\em same} radius as the original balls. Indeed, if one
allows a doubling of the radius, the corresponding inequality is
trivial, but insufficient for taking the thermodynamic limit in the
entropy. Thus, our bound is an essential tool for showing the
existence of infinite volume limits in topologies where the
``matching'' of functions needs some care.

Once all these tools are in place, we can easily repeat the proofs of
the existence of the topological entropy 
using the methods developed in [CE1] and [CE2].

The paper is structured as follows. In \sec{coercive} we bound the
flow in time, using localized versions of coercive functionals as
introduced by Feireisl [F]. The main result is \clm{f1} and its
corollary
\equ{good3} and \equ{good4} which show that the solution to
\equ{11a}
is well behaved in $\HH_{\delta,\loc}^3$. In \sec{linearized} and
\sec{prep} we study the linear part of \equ{problem} {\em localized in
coordinate and momentum space}. We next study the decay of the high
momentum part in \sec{high} and \sec{highfreq}. This allows, in
\sec{full} to study the time evolution of differences of two solutions
of \equ{problem}, in other words, we control now the {\em deformation}
of balls (in the topology of $\HH_{\delta,\loc}^2$).
In \sec{bounds} (\clm{finite} and \equ{balls2})
we show how to cover the attracting set $\GG$ with balls as explained
in \equ{bbbb} above. In \sec{subadditivity} we deal with the
technically delicate submultiplicativity bound mentioned before.
Finally, in \sec{ee} and \sec{top} we show without effort the bound
\clm{ee} on the $\epsilon $-entropy (per unit length) and the
\clm{exist}
which shows the existence of the topological entropy per unit length.

\SECT{coercive}{Coercive functionals}

In this section, we study some functionals which control the flow in
time. The first part of this material is an adaptation from the work of
Feireisl[F].
We consider here the problem \equ{start1} in the form
$$
\eqalign{
\dot u \,&=\,v~,\cr
\eta^2\dot v \,&=\,-v +u'' +U'(u)~,\cr
}
\EQ{problem}
$$
where we take $U(s)=\frac{s^2}2-\frac{s^4}4$, but many other choices
are of course possible.
To simplify things, we {\em assume throughout that $0<\eta<1$}, and in
fact, in subsequent sections we will assume $\eta<\eta_0$ for some
small $\eta_0$.
We shall use throughout a localization function $\hh_\epsilonp $ which
depends on a small 
parameter $\epsilonp $, to be determined later on. The constant
$\epsilonp$ will 
only depend on the coefficients of \equ{problem} (but not on
$\eta<1)$.

We set
$$
\hh_{\epsilonp }(x)\,=\,{1\over (1+\epsilonp ^2x^2)^2}~.
\EQ{hh}
$$
Note that $\epsilonp \int\dx\hh_{\epsilonp}$ is independent of
$\epsilonp$.
\REMARK We will only use values of $0<\epsilonp \le\HALF$ and this will be
tacitly assumed in all the estimates.

Using $\hh_{\epsilonp }$, we introduce the norm 
$$
\|(u,v)\|^2_{\hh_{\epsilonp},1}\,=\,
\int\dx \hh_{\epsilonp}(x) \bigg(\eta^2 v^2+u^2+(u')^2
\bigg )(x)~.
\EQ{gdd}
$$
We also need a translation invariant topology on $(u,v)$.
Let $\hh_{\epsilonp,\xi}(x)=\hh_{\epsilonp}(x-\xi)$.
\CLAIM{Definition}{bbloc1}We define the norm
$$
\|(u,v)\|_{\alpha
,\loc,1}\,=\,\sup_{\xi\in\real}\|(u,v)\|_{\hh_{\epsilonp,\xi},1}~,
\EQ{alphaloc1} 
$$
and the space $\HH_{\epsilonp,\loc}^1$ is defined by
$$
\HH_{\epsilonp,\loc}^1\,=\,\bigg\{ (u,v) ~\bigg |~\|(u,v)\|_{\alpha
,\loc,1}
<\infty 
\bigg\}~.
\EQ{bb1def}
$$

The norm introduced above is not very convenient for estimates, and
thus we introduce as in Feireisl[F] 
the quantity $F_0$ (which is not a norm) by
$$
F_0(u,v)\,=\,\alpha \int \dx \hh_{\epsilonp}(x)\biggl(\eta^2 v^2(x) + (u'(x))^2 +
V(u(x))+\eta^2 \nut1 u(x) v(x)\biggr)~.
\EQ{F00}
$$
Here,
$V$ is chosen 
such that
$$
U'(x)+2V'(x)-\nut1 \eta^2 x\,=\,0~,
\EQ{G0}
$$
with $V(0)=0$. 
Note that $U(x)\to-\infty $ as $|x|\to\infty$ at a rate $\OO(x^4)$,
and therefore $V(x)\to+\infty $ at a rate $\OO(x^4)$. 
In particular,
$$
V(x)\,\ge\,x^2\text{ for sufficiently large }|x|~.
\EQ{gbound}
$$
The following bound can be found in Feireisl[F].
\CLAIM{Lemma}{E0}There are constants $a_0>0$ and $b_0>0$ (independent
of $\eta$ for $0<\eta<1$)
for which one has the inequality
$$
\partial _t F_0(u_t,v_t)
\,\le\,
-a_0F_0(u_t,v_t)+b_0~,
\EQ{dtE00}
$$
where $u_t(x)=u(x,t)$, $v_t(x)=v(x,t)$ is the solution of
\equ{problem}.
 
This bound can be used to bound $\|(u,v)\|_{\hh_{\epsilonp},1}$. 
Recall that $V$ diverges like $|x|^4$. Using the bound
$$
\eta^2\nut1 |uv|\,\le\, {\eta^2 v^2\over 2} + {\eta^2 \nut1\nut1 u^2\over
2}~,
$$
this implies
$$
V(u)+\eta^2\nut1 uv \,\ge\, -{\eta^2v^2\over 2} +u^2 -\C{g}~,
$$
for some constant $\CX{g}$.
Therefore, one has the inequality 
$$
\|(u,v)\|_{\hh_{\epsilonp},1}^2\,\le\,2 F_0(u,v) +\C{g1}~.
$$
Using \equ{dtE00}
we conclude
\CLAIM{Lemma}{E00}There is a constant $\C{g2}$ (independent of
$0<\eta<1$) for which the following holds.
Assume that $F_0(u_0,v_0)<\infty $.
Then, for all $t>0$ one has $\|(u_t,v_t)\|_{\hh_{\epsilonp},1}<\infty $
and there is a $T=T(u_0,v_0)$ for which the solution $(u_t,v_t)$
of \equ{problem}
with initial data $(u_0,v_0)$ satisfies for all $t>T$:
$$
\|(u_t,v_t)\|_{\hh_{\epsilonp},1} \,\le\, \CX{g2}~.
\EQ{dtE0}
$$

We can extend this result to the topology of $\HH_{\alpha ,\loc}^1$.
Let $u_{0,\xi}(x)=u_0(x-\xi)$ and similarly for $v_0$.
\CLAIM{Proposition}{linfty}There is a constant $\C{h}$ (independent of
$0<\eta<1$) for which the following holds.
Assume that $\sup_{\xi\in\real} F_0(u_{0,\xi},v_{0,\xi})<\infty $.
Then there is a $T=T(u_0,v_0)$ for which the solution $(u_t,v_t)$
of \equ{problem}
with initial data $(u_0,v_0)$ satisfies for all $t>T$:
$$
\|(u_t,v_t)\|_{\alpha ,\loc,1}\,\le\,\CX{g2}~,
\EQ{thisisit}
$$
and
$$
\|u_t\|_\infty \,\le\,\CX{h}~.
$$
\PROOF Consider the quantities $F_{0,\xi}$ defined by
replacing $\hh_{\epsilonp}(x)$ by its translate
$\hh_{\epsilonp}(x+\xi)$ in Eq.\equ{F00}. Then,
$F_0(u_{t,\xi},v_{t,\xi})=
F_{0,\xi}(u_t,v_t)$.
Clearly, for every $\xi$ we have
$$
\partial _t F_{0,\xi}(u_t,v_t)\,\le\,-a_0 F_{0,\xi}(u_t,v_t) +b_0~,
$$
since \equ{problem} does not depend explicitly on $x$.

It follows from the above that if
$\sup_\xi F(u_{0,\xi},v_{0,\xi})<\infty $
there is a 
finite time $T$ after which 
$$
\|(u_t,v_t)\|_{\alpha ,\loc,1}\,\le\,\CX{g2}~.
\EQ{aux}
$$
This proves \equ{thisisit}.  
To conclude that $u$ is bounded we need the following easy
\CLAIM{Lemma}{Sobolev}There is a constant
$\C{Sobolev}=\CX{Sobolev}(\delta)$ such that  
$$
\sup_{x\in[-1,1]}|f(x)|\,\le\, \CX{Sobolev} \|f\|_{\hdd,1}~.
\EQ{Sobolev}
$$

\PROOF From the explicit form of $h_\delta $ we conclude
that
$$
\int_{-1}^1\dx|f(x)|^2\,\le\,(1+\delta ^2)^2 \int \dx \hdd |f(x)|^2~,
$$
and similarly
$$
\int_{-1}^1\dx|f'(x)|^2\,\le\,(1+\delta ^2)^2 \int \dx \hdd |f'(x)|^2~.
$$
The result follows from the standard Sobolev inequality.
The proof of \clm{Sobolev} is complete.\QED
Using this lemma, and observing that the $\|\cdot\|_{\alpha ,\loc,1}$
norm is translation invariant we conclude immediately from \equ{aux}
that there is
a constant $\CX{h}$ 
for which
$$
\sup_x |u(x,t)|\,\le\, \CX{h} ~,
\EQ{ubound}
$$
for all $t>T$. The proof of \clm{linfty} is complete.\QED

We next deal with the slightly more complicated bounds on the
{\em spatial derivatives} of $u$ and $v$. Let $w=u'$ and let $z=v'$. They
satisfy the equations
$$
\eqalign{
\dot w_t \,&=\,z_t~,\cr
\eta^2\dot z_t \,&=\,-z_t +w_t'' +U''(u_t)w_t~,\cr
}
\EQ{wz}
$$
where $U''(s)=1-3s^2$, for the $U$ we have taken above. 
We consider initial data $(w_0,z_0)$ which will be bounded later and 
assume (in view of \clm{linfty}) that
$\|(u_0,v_0)\|_{\alpha ,\loc,1}\le\CX{g2}$ which implies $\sup_x|u(x,t)|\le\CX{h}$ for all $t>0$. 

We are going to bound the growth of $(w,z)$ as a function of time.
We introduce a positive constant
$\mu$ (which we fix later) and set
$$
F_1(w,z)\,=\,
\alpha \int \dx \hh_{\epsilonp}(x)\biggl(\eta^2 z^2(x) + (w'(x))^2 +
\eta^2 \mu  w(x) z(x)\biggr)~.
\EQ{F10}
$$
When no confusion is possible, we henceforth write $\int f$ for
$\int\dx f(x)$.
One finds for $(w_t,z_t)$ satisfying \equ{wz}: 
$$
\eqalign{
\HALF\partial _t F_1(w_t,z_t)\,&=\,
\alpha \idx  \hh_{\epsilonp}\bigg(
\eta^2 z_t \dot z_t + w_t' \dot w_t' +{\eta^2\mu  \over 2} \dot w_tz_t+{\eta^2\mu 
\over 2}w_t\dot z_t
\bigg)\cr\,&=\,
\alpha \idx  \hh_{\epsilonp}\bigg(
-z_t^2+z_tw_t''+U''(u_t)z_tw_t + w_t'z_t' +{\eta^2\mu  \over 2} z_t^2\cr
&~~~~~-{\mu 
\over 2}w_tz_t +{\mu 
\over 2}w_tw_t''+{\mu 
\over 2}U''(u_t)w_t^2
\bigg)\cr
\,&=\,
\alpha \idx  \hh_{\epsilonp}\bigg(
-z_t^2(1-{\eta^2\mu \over 2})+U''(u_t)z_tw_t \cr
&~~~~~-{\mu 
\over 2}w_tz_t -{\mu 
\over 2}(w_t')^2+{\mu 
\over 2}U''(u_t)w_t^2
\bigg)\cr
&~~~~~-\alpha \idx  \hh_{\epsilonp}'\bigl(z_tw_t'+{\mu \over 2}w_tw_t'\bigr)~.\cr
}
\EQ{long}
$$
Note now that by the definition \equ{hh} of
$\hh_{\epsilonp }$ and the restriction $\epsilonp\le \HALF$, we find
that the quotient $|\hh_{\epsilonp}'/\hh_{\epsilonp}|$ is bounded by
$2\epsilonp\le 1$.
Using this, we get
$$
|\hh_{\epsilonp}' z_t w_t'|\,\le\,\epsilonp \hh_{\epsilonp}\bigl(z_t^2+(w_t')^2\bigr )~,
$$
and
$$
|\mu \hh_{\epsilonp}'w_t w_t'|\,\le\,\epsilonp \mu \hh_{\epsilonp}\bigl( w_t^2 +(w_t')^2\bigr )~.
$$
We will also use the bound $|w_tz_t|\le \epsilonp ^{-1} w_t^2+\epsilonp z_t^2$.
Finally note that there is a constant $\C{f2}$ for which
$$
\sup_{|s|\le \CX{h} } |U''(s)| \,\le\,\CX{f2}~.
$$
Combining these bounds with the last equality of \equ{long}, we
get for times $t>0$:
$$
\eqalign{
\HALF\partial _t F_1(w_t,z_t)\,&\le\,
-\alpha \idx  \hh_{\epsilonp} 
z_t^2\bigl(1-{\eta^2\mu \over 2}-\epsilonp -\mu /2-\epsilonp  \CX{f2}\bigr)\cr
&~~~~
-\alpha \idx  \hh_{\epsilonp} (w_t')^2\bigl(\mu /2-\epsilonp \mu /2 -\epsilonp \bigr)\cr
&~~~~+\bigl(\epsilonp \mu +\mu /2+\mu \CX{f2}/2+\epsilonp ^{-1}\CX{f2}\bigr)
\alpha \idx  \hh_{\epsilonp} w_t^2
~. 
\cr
}
\EQ{long2}
$$
It is clear that if we choose $\epsilonp $ and $\mu $ sufficiently
small
(but independent of $\eta$ for small enough $\eta$),
then we get, for some (large) constant $\C1$,
$$
\eqalign{
\HALF\partial _t F_1(w_t,z_t)\,&\le\,
\alpha \idx  \hh_{\epsilonp}  \bigg(
-{\mu\over 4}\bigl ( z_t^2 + (w_t')^2\bigr )
+\CX1 w_t^2\bigg )~.\cr
}
\EQ{long3}
$$
Note now that for small $\mu>0$ and $\eta>0$,
$$
\eta^2\mu wz \,\le\,(1-\eta^2) z^2+{\eta^2\mu\over 4(1-\eta^2)} w^2~,
$$
which is equivalent to
$$
-z^2\,\le\,-\eta^2z^2 -\eta^2\mu wz +{\eta^2\mu\over 4(1-\eta^2)}
w^2~.
$$
Therefore, \equ{long3} leads to
$$
\eqalign{
\HALF\partial _t F_1(w_t,z_t)\,&\le\,
-{\alpha\mu\over 4} \idx  \hh_{\epsilonp}  
\bigl ( \eta^2 z_t^2 + (w_t')^2+\eta^2 \mu w_t z_t\bigr )
+\C{777}  \idx  \hh_{\epsilonp}w_t^2~.\cr
}
$$
The last term is bounded by \equ{thisisit}, and therefore we find:
$$
\partial _t F_1(w_t,z_t)\,\le\,-a_1 F_1(w_t,z_t) + b_1~,
$$
for some finite positive $a_1$ and $b_1$.
Using again the methods leading to \clm{linfty}, we obtain
\CLAIM{Theorem}{f1}There are constants $\C{7777}$, $\C{77}$ and $\C7$
(independent of 
$\eta<1$) for which the following holds.
Assume $\sup _{\xi\in\real}F_{0}(u_{0,\xi},v_{0,\xi})<\CX{7777}  $ and
$\sup _{\xi\in\real}F_{1}(u_{0,\xi},v_{0,\xi})<\CX{7777}  $.
Then the solution $(u_t,v_t)$
of \equ{problem}
with initial data $(u_0,v_0)$ satisfies for all $t>0$:
$$
\|(u_t,v_t)\|_{\alpha ,\loc,2}\,\le\,\CX{77}~,
\EQ{good0} 
$$
and
$$
\sup _{x\in\real}\bigl(|u(x,t)| + | u'(x,t)| +|v(x,t)|\bigr)\,\le\,\CX7~.
\EQ{good}
$$

\REMARK The technique used above can be extended to show that any
derivative of $u(x,t)$ and $v(x,t)$ is eventually bounded (if the
potential $U$ is sufficiently differentiable and the initial data are
sufficiently regular). The details are left to
the reader. We will in fact use bounds on the second derivative at
some later point in the argument, {\it i.e.},
bounds of the form
$$
\|(u_t,v_t)\|_{\alpha ,\loc,3}\,\le\,\CX{77}~,
\EQ{good3} 
$$
and
$$
\sup _{x\in\real}\bigl(|u(x,t)| + | u'(x,t)| + | u''(x,t)| +|v(x,t)| +|v'(x,t)|\bigr)\,\le\,\CX7~.
\EQ{good4}
$$

\SECT{linearized}{The linearized evolution}

In this section, we study the linear part of the solution.
By this we mean solutions of the equation
$$
\eqalign{
\dot \gu \,&=\,\gv ~,\cr
\eta^2\dot \gv \,&=\,-\gv +\gu '' ~.
\cr
}
\EQ{diff1}
$$
It will be
useful to rewrite this system of equations as
$$
\V{\dot \gu}{\dot \gv}\,=\,\LL\V{\gu}{\gv}~.
\EQ{llform}
$$

Next, we go through a second round of estimates, similar to the ones in
\sec{coercive}, 
to see how fast  $(\gu ,\gv )$ can grow.
We use again the cutoff function
$$
\hhh_\delta(x)\,=\,{1\over (1+\delta^2x^2)^2}~,
$$
with a parameter $\delta $ different from $\epsilonp$
which will be fixed in \sec{highfreq}.
We consider the functional
$$
H\,=\,\idx  \hhh_\delta\bigl(\eta^2 \gv ^2  +\gu ^2+(\gu ')^2\bigr) ~\,=\,\|(\gu,\gv)\|_{\hhh_\delta,1}~, 
\EQ{hfirst}
$$
and proceed to bound it.
One gets for the solution $(\gu_t,\gv_t)$ of \equ{llform} : 
$$
\eqalign{
\HALF \partial _t  H(\gu_t,\gv_t)\,=\,
\int &\dx \hhh_\delta\bigl(\eta^2\gv_t \dot \gv_t  + \gu_t \dot \gu_t +\gu_t '\dot \gu_t ' \bigr)\cr
\,=\,\int &\dx \hhh_\delta\bigl(- \gv_t ^2 +\gv_t \gu_t ''+\gv_t \gu_t +\gu_t '\gv_t '\bigr) \cr
\,=\,\int &\dx \hhh_\delta\bigl(-\gv_t ^2+\gv_t \gu_t \bigr) -\idx \hhh_\delta' \gv_t  \gu_t'~.\cr
}
\EQ{verylong}
$$
Observe that by construction,
$$
|\hhh_\delta'/\hhh_\delta|\,\le\,1~,
$$ 
when $\delta <\HALF$ (which we assume throughout).
We use now (since $\gu$ and $\gv$ are real):
$$
|\gv_t |\bigl(|\gu_t |+|\hhh_\delta'\gu_t '/\hhh_\delta|\bigr)
\,\le\, \HALF \gv_t ^2 +\gu_t ^2 +(\gu_t ')^2~.
$$
We can use this inequality to bound the mixed terms in the last
equality of \equ{verylong}, and compensate the term
$\HALF(\gv_t ^2)$ with the negative term $-\gv_t ^2$.
Therefore, we have shown that there is a constant $\C{H}$ {\em
independent of $\eta<1$} for which
$$
\partial _t H(\gu_t ,\gv_t )\,\le\,\CX{H} H(\gu_t ,\gv_t )~.
\EQ{Hfinal}
$$

We next define a space in which both the functions and their derivatives are
controlled. 
This is the space in which our final bounds will be spelled out.
Using still the cutoff function
$$
\hdd(x)\,=\,{1\over (1+\delta^2x^2)^2}~,
$$
we define the norm
$$
\|(u,v)\|^2_{h_{\delta},2}\,=\,\idx  \hdd \bigg(\eta^2 \bigl(v^2+(v')^2\bigr)+u^2+2 (u')^2
+(u'')^2\bigg ) ~.
\EQ{ldd}
$$
\CLAIM{Definition}{H2}The Hilbert space $\HH_{\hdd}^2$ is defined by
$$
\HH_{\hdd}^2\,=\,
\bigg\{ (u,v) ~\bigg|~  \|(u,v)\|_{h_{\delta},2}<\infty
\bigg\}~.
\EQ{h2def}
$$

We also need the translates of $\hdd$ to define
a translation invariant topology.
Let $h_{\delta,\xi}(x)=\hdd(x-\xi)$.
\CLAIM{Definition}{bbloc2}The space $\HH_{\delta,\loc}^2$ is defined by
$$
\HH_{\delta,\loc}^2\,=\,\bigg\{ (u,v) ~\bigg|~  \sup_\xi\|(u,v)\|_{h_{\delta,\xi},2}<\infty
\bigg\}~.
\EQ{bb2def}
$$

\REMARK Note that 
$$
\|(u,v)\|_{\hdd,2}^2\,=\,\|(u,v)\|_{\hdd,1}^2+
\|(u',v')\|_{\hdd,1}^2~.
$$

Next, we observe that $\LL$ and $\partial _x$ commute. Therefore, the
bounds on $\|(\gu_t,\gv_t)\|_{h_{\delta},1}$ can be extended
immediately to similar bounds on
$\|(\gu_t,\gv_t)\|_{h_{\delta},2}$
and we get from \equ{Hfinal}: 
\CLAIM{Theorem}{frakuv}There is a constant $\CX{H}$ such that if 
$(\gu_t,\gv_t)$ solves \equ{diff1} then 
$$
\|(\gu_t,\gv_t)\|_{\hdd,2}^2\,\le\,e^{\CX{H}t}\|(\gu_0,\gv_0)\|_{\hdd,2}^2~.
\EQ{growth}
$$

\SECT{prep}{Momentum localization}

Let $m$ be a bounded positive function of $x$ which decays faster than
any inverse power as $|x|\to \infty $.
We define
the convolution operator $M_a$ by
$$
\bigl(M_af\bigr)(x)\,=\,\int \d{y} m\bigl(a(x-y)\bigr) f(y)~.
\EQ{Qdef}
$$
Let again
$$
\hdd(x)\,=\,{1\over (1+\delta^2x^2)^2}
~.
\EQ{hdef}
$$
\CLAIM{Lemma}{q}There is a constant $\C{M}$ such that
if $\delta>0$ and $a>0$,
then the operator $M_{a}$ is bounded on $\L^2(\hdd\,\,\dx)$,
with norm bounded by
$$
\| M_a\| _\delta\,\le\,\CX{M}{1+{\delta^2\over a^2}\over a}~.
\EQ{boundona}
$$

\PROOF We will prove the result by bounding the operator $\hat M_a$
with integral kernel
$$
\hdd^{1/2}(x)\,m\bigl(a(x-y)\bigr)\,\hdd^{-1/2}(y)
\EQ{kernel}
$$
on $\L^2({\rm d}x)$.
Writing $1=\chi(2|x|>|y|)+\chi(2|x|<|y|)$, with characteristic
functions $\chi$, and multiplying the kernel \equ{kernel} with
them, we induce
a decomposition of this operator as a sum of two pieces.

Fix $|x|$. We consider first the integration region
$|y|<2|x|$. In that region, we have a bound
$$
0\,<
\,{\hdd^{1/ 2}(x)\over \hdd^{{1/2}}(y)}\,\le\,4~.
$$
Therefore, in this region, the integral kernel is bounded by $4
|m(a(x-y))|$. Since
$m $ decreases like an arbitrary power we get for every
$\ell>0$
a bound
$$
\int_{|y|<2|x|}\kern-1.5em\d{y} { h_\delta^{1/2} (x) |m\bigl(a
(x-y)\bigr )| \over
h_\delta(y)^{1/2}}|f(y)|\,\le\,
4K_\ell\int_{|y|<2|x|}\kern-1.em\d{y} {1\over
(1+a|x-y|)^\ell}|f(y)|~.
$$
When $\ell>1$, then $a/(1+a|x|)^\ell$ is bounded in $\L^1$ (uniformly in
$a$) and therefore,
Young's inequality shows that this piece of $\hat M_af$ is bounded in
$\L^2({\rm d}x)$,
with norm less that $\C{Q<}a^{-1}\|f\|_2$, and $\CX{Q<}$ independent of
$\delta$ and $a$.
In the region $|y|\ge2|x|$, we use $|x-y|>|y|/2$ and $|x-y|>|x|$.
Therefore, using $\hdd(x)\le1$, we find for $\ell>3$,
$$
\int_{|y|\ge2|x|}\kern-1.5em\d{y} { h_\delta^{1/2} (x)
|m\bigl(a(x-y)\bigr )| \over
h_\delta(y)^{1/2}}|f(y)|
\,\le\,
2^{(\ell+2)\over2}K_\ell \int_{|y|\ge2|x|}\kern-1.em\d{y} {1+\delta^2y^2\over
(1+a|x|)^{\ell -2\over 2}(1+a|y|)^{\ell +2\over 2}}|f(y)|~.
$$
Note that 
$$
{1+\delta^2y^2\over (1+a|y|)^2}\,\le\,1+{\delta^2\over a^2}~.
$$
Using the Schwarz inequality yields a bound $(1+\delta^2/a^2)/a$
on the second piece of
$\|M_a\|_\delta$. Combining the two pieces completes
the proof of \clm{q}.\QED

We need later the following variant of this result:
Let $\theta=\theta(k)$ be a smooth characteristic function which equals 1
for $|k|\le1$ and 0 for $|k|>2$. Let $q_{a}(x)=a\tilde
\theta(ax)$, with $\tilde \theta$ the inverse Fourier transform of
$\theta$.
For $a>0$, let $\QQ_{a}$ be the convolution operator defined by
$$
\bigl(\QQ_{a}f\bigr)(x)\,=\,\int\d{y} q_{a}(x-y) f(y)~.
\EQ{Qdef0}
$$
This operator is a substitute for a projection onto momenta less than $a$.
Setting $m(x)=\tilde \theta(x)$, we get from \clm{q}:
\CLAIM{Corollary}{Qa}There is a constant $\CX{M}$ such that if
$\delta>0$ and $a>0$ then the operator $\QQ_a$ is bounded on
 $\L^2(\hdd\,\,\dx)$,
with norm bounded by
$$
\| \QQ_a\|_\delta \,\le\,\CX{M}\bigl({1+{\delta^2\over a^2}}\bigr )~.
\EQ{boundonq}
$$

\SECT{high}{High momentum bounds}

We consider again the function $\hdd$ as defined in \equ{hdef}, and we
study functions $u$ for which $\int \dx\hdd(x) \bigl(|u(x)|^2 +
|u'(x)|^2\bigr)<\infty$.
Such functions have a Fourier
transform $\tilde u$ in the sense of tempered distributions, and we define now 
$$
\KK_a\,=\,\bigg\{ u  ~ \bigg |~ \int \dx\hdd(x) \bigl(|u(x)|^2 +
|u'(x)|^2\bigr)<\infty \text{ and } \supp \tilde u \in \real
\setminus (-a,a)\bigg\}~.
$$
Thus, apart from not being defined as a function, the Fourier
transform $\tilde u$ of a $u\in \KK_a$ has support at momenta larger than $a$.
If $\hdd(x)\equiv 1$ and $u'\in \L^2({\rm d}x)$, then, obviously, for
$u\in \KK_a$,
one has $\idx  |u|^2\,\le\,a^{-2}\idx  |u'|^2$.
The following proposition whose elegant proof was kindly provided by
{\em H. Epstein}, shows that the cutoff function $\hdd$ does
only
moderately change this property.
\CLAIM{Proposition}{kaa}Assume that $a>0$ and
$\delta>0$.
There is a constant $\nu(a,\delta)<\infty$ such that
for all $u\in\KK_a$ one has the inequality
$$
\int \dx \hdd (x) |u(x)|^2 \,\le\,
\nu (a,\delta)\int \dx \hdd (x) |u'(x)|^2~.
\EQ{kaa}
$$
There is a constant $\C{nu0}>0$ such that one can choose
$$
\nu(a,\delta)\,=\,\CX{nu0}{(1+{\delta^2\over a^2})^2\over a^2}~.
\EQ{numore}
$$

\REMARK We will need the result only for $\delta<a$, so that we can
use the simpler bound
$$
\nu(a,\delta)\,\le\,{\C{nu}\over a^2}~.
\EQ{nufinal}
$$
\PROOF Let $\theta$ be a smooth characteristic function which equals 1
for $|k|\le\HALF$ and 0 for $|k|>1$. Let $u\in\KK_a$. Since the distribution
$\tilde u$ has support in the complement of the interval $(-a,a)$, and
the Fourier transform $\widetilde {u'}$ of the derivative is $ik\tilde
u(k)$, we see that
$$
\tilde u(k)\,=\,{1-\theta(k/a)\over ik}\,\,\widetilde {u'}(k)~.
$$
Define next 
$$
\tilde m(k)\,\equiv\,{1-\theta(k)\over ik}~.
$$
The (inverse) Fourier transform, $m$, of $\tilde m$ decreases faster
than any power of $|x|$ at infinity.
If we let $m_a(x)=m(ax)$, then $\tilde m_a(k)=\tilde m(k/a)/a.$
Thus, it follows with the notation of \sec{prep} that
$$
u(x)\,=\,\bigl(M_a(u')\bigr)(x)~.
$$
By \clm{q}, we conclude that
$$
\|u\|_\delta^2\,=\,\idx  \hdd |u|^2 \,=\,\idx  \hdd
|M_a(u')|^2\,=\,\|M_a (u')\|_\delta^2
\,\le\,\CX{M}^2 {(1+{\delta^2\over a^2})^2\over a^2}\idx  \hdd |u'|^2~,
$$
and the claim \equ{numore} follows.\QED

\SECT{highfreq}{The linear high frequency part}

We begin by defining the projection onto high frequencies, on a space
with weight $\hdd(x)=(1+\delta^2x^2)^{-2}$.
We first recall the notion of projection onto low frequencies from \sec{prep}.
Denote by $\widetilde\theta$  a smooth characteristic function, equal to 1 for $|k|\le1$
and vanishing for $|k|>2$. We fix now a (large) cutoff scale $k_*$ and
we define as before
$$
q_{k_*}(x)\,=\,k_*  \theta(k_* x)~,
$$
and
$$
\bigl(\QQ_{k_*}f\bigr)(x)\,=\,\int \d{y} q_{k_*}(x-y) f(y)~.
\EQ{Qdef2}
$$
In \clm{Qa}, we showed that on $\L^2(\hdd\,\,\dx)$, the operator
$\QQ_{k_*}$ is bounded by $\CX{M}\bigl({1+{\delta^2\over k_*^2}}\bigr
)$.
Therefore, the projection onto high momenta
$$
\PP_{k_*}\,=\,1-\QQ_{k_*}~,
\EQ{pdef}  
$$
is also bounded on that space.
Henceforth, we shall assume $\delta<k_*$, and thus we get immediately the bound
\CLAIM{Lemma}{q2}There is a constant $\C{Q}$ such that
if $k_*>\delta>0$,
then the operator $\PP_{k_*}$ satisfies
$$
\|\PP_{k_*}\|_\delta\,\le\,\CX{Q}~,
\EQ{Pbound}
$$
as a map on  $\L^2(\hdd\,\,\dx)$.

\CLAIM{Lemma}{QonB}There is a constant $\C{QQ}$ such that
for $k_*>\delta>0$
the operator $\PP_{k_*}\oplus \PP_{k_*}$ is bounded in norm by $\CX{QQ}$
as a map from $\HH_{\delta,\loc}^2$ to itself.

\PROOF We have already checked in \clm{q2} that $\PP_{k_*}$ is bounded
on $\L^2(\hdd\,\,\dx)$. Note that $\PP_{k_*}$ is a convolution operator
and so
$\PP_{k_*}$ and $\partial_x$ commute, and the extension of the
result to $\HH_{\delta,\loc}^2$ (as defined in \clm{bbloc2}) follows at once.\QED

So far, we have argued that $\PP_{k_*}$ is bounded.
We will now use 
the high momentum bound of \sec{high} with
$a=k_*$, and $k_*\le\eta^{-1}$ to show that the semi-group generated
by the free evolution (see below) is a (strong) contraction.
In fact, we will show that the contraction rate is
$$
\displaylines{
\OO(k_*^2)\text{ as long as } k_* <\eta^{-1}~,\cr 
\OO(\eta^{-2})\text{ for any cutoff }k_*\ge\eta^{-1}~.
\cr
}
$$
This behavior is typical of the mixed
parabolic-hyperbolic problems we consider here, since the linearized
evolution, written in momentum space, has the generator
$$
\left(
\matrix{0& 1\cr -\eta^{-2}k^2 &-\eta^{-2}}
\right )
$$
with eigenvalues 
$$
\lambda_\pm\,=\, {-1\pm\sqrt{1-4k^2\eta^{-2}}\over
2\eta^{2}}~. 
$$
One can see from the expression for the eigenvalues that the real part
never gets more negative than $-\OO(\eta^{-2})$.
We study now the properties of the operator $\LL$ defined as in
Eqs.\equ{diff1} and \equ{llform} by
$$
\LL\V{\gu }{\gv }\,=\,\V{\gv }{\eta^{-2}\bigl(-\gv +\gu ''\bigr)}~.
$$

We introduce parameters $\gamma >0$, and $\delta>0$ which will be
fixed later and
we consider the functional $J$:
$$
J\,=\,J_{\hdd ,\gamma }(\gu ,\gv )\,=\,\|(\gu ,\gv )\|_{\hdd,2}^2+\eta^2\gamma \int\dx \hdd (x) \bigl(\gu (x) \gv (x) +
\gu '(x)\gv '(x)\bigr )~,
\EQ{jdef0}
$$
where the norm $\|(\gu ,\gv )\|_{\hdd,2}$ was defined in Eq.\equ{ldd}:
$$
\|(\gu ,\gv )\|^2_{h_{\delta},2}  \,=\,\idx  \hdd \bigg(\eta^2 \bigl(\gv ^2+(\gv ')^2\bigr)+\gu ^2+2 (\gu ')^2
+(\gu '')^2\bigg )~.
$$
Consider the solution $(\gu_t,\gv_t)$ of \equ{llform}: 
Then, writing $J$ for $J_{\hdd, \gamma }$, we find
$$
\eqalign{
\HALF \partial_t J(\gu_t,\gv_t) \,&=\,
\idx  \hdd  \bigg(
\eta^2\bigl( \gv_t \dot \gv_t +\gv_t '\dot \gv_t '\bigr )+ \gu_t \dot \gu_t  +2\gu_t '\dot \gu_t ' + \gu_t ''\dot \gu_t ''
\biggr)\cr
&+\HALF\eta^2\gamma \idx  \hdd  \bigl(\gu_t \dot \gv_t +\gv_t \dot \gu_t +\gu_t '\dot
\gv_t '+\gv_t '\dot \gu_t '\bigr)\cr
\,&=\,
\idx  \hdd  \bigg(
-\gv_t ^2 +\gv_t \gu_t '' -(\gv_t ')^2 +\gv_t '\gu_t ''' +\gu_t \gv_t  +2\gu_t '\gv_t ' +\gu_t ''\gv_t ''
\bigg)\cr
&+\HALF\gamma  \idx 
\hdd \bigl(-\gu_t \gv_t +\gu_t \gu_t ''+\eta^2\gv_t ^2-\gu_t '\gv_t '+\gu_t '\gu_t '''+\eta^2(\gv_t ')^2\bigr)\cr
\,&=\,
\idx  \hdd  \bigg(
-\bigl(\gv_t ^2+(\gv_t ')^2\bigr)(1-\HALF\eta^2\gamma )
+\bigl(\gu_t \gv_t +\gu_t '\gv_t '\bigr)(1-\HALF\gamma )\bigg)\cr   
&-\HALF\gamma \idx  \hdd  \bigl((\gu_t ')^2+(\gu_t '')^2\bigr)-\idx 
\hdd '\bigl(\gv_t \gu_t '+\gv_t '\gu_t ''+\HALF\gamma \bigl(\gu_t \gu_t '+\gu_t '\gu_t ''\bigr)\bigr)  
~.\cr
}\EQ{hdd}
$$
By construction, we have $|\hdd'/\hdd|\le \C{dd}\delta$, and therefore
the last integral in \equ{hdd} can be bounded (in modulus) by
$$
\CX{dd}\delta\idx  \hdd
\bigg(\gv_t ^2+(\gu_t ')^2+(\gv_t ')^2+(\gu_t '')^2+\HALF\gamma \bigl(\gu_t ^2+2(\gu_t ')^2+(\gu_t '')^2\bigr)
\bigg)~. 
$$
Thus, we find
$$
\eqalign{
\HALF\partial _t J(\gu_t,\gv_t)\,\le\,
&-\idx  \hdd \bigl(\gv_t ^2+(\gv_t ')^2\bigr
)\bigl(1-\HALF\eta^2\gamma -\CX{dd}\delta\bigr)\cr
&-\HALF\gamma \idx  \hdd
\bigl((\gu_t ')^2+(\gu_t '')^2\bigr)\bigl(1-{\CX{dd}\delta\over
\gamma }-2\CX{dd}\delta \bigr)\cr
&+\HALF\CX{dd}\delta\gamma \idx  \hdd \gu_t ^2\cr
&+\idx  \hdd\bigl(\gu_t \gv_t +\gu_t '\gv_t '\bigr)(1-\HALF\gamma ) 
~.\cr
}\EQ{hdd2}
$$
Recall that $\eta>0$ is given, and that we want to prove results for all
$\eta<\eta_0$, where $\eta_0$ is our (only) small parameter.

We rewrite the last integral in \equ{hdd2} as
$$
-\EIGHT\eta^2\gamma ^2\idx  \hdd\bigl(\gu_t \gv_t +\gu_t '\gv_t '\bigr)
+\idx  \hdd\bigl(\gu_t \gv_t +\gu_t '\gv_t '\bigr)(1+\EIGHT\eta^2\gamma ^2-\HALF\gamma )~.
\EQ{split}
$$ 
We introduce now the first restrictions on $\eta$ and $k_*$:
Fix
$$
\eta_0\,\le\,{1\over \sqrt{40}}~,
\EQ{eta0}
$$
and
$$
 k_0\,\ge\,\sqrt{40\CX{nu}}~.
\EQ{defk*}
$$
These bounds will be made more stringent below.
We shall always require
$$
0\,<\,\eta\,<\,\eta_0~,\quad \text{and} \quad
\infty \,>\,k_*\,>\, k_0~.
$$
We next define
$$
\eqalign{
\gamma \,&\equiv\, \min(\eta^{-2},k_*^2/\CX{nu})/320~,\cr
}
\EQ{ddx0} 
$$
and we choose the space-cutoff parameter $\delta $ as
$$
\delta\,=\,\min\bigl (1/2,1/(40\CX{dd})\bigr )~.
\EQ{deltadef}
$$

Note that $\gamma$ is essentially the inverse of the dispersion
law as explained at the beginning of this section.

With the above requirements we find $\gamma >2$ and
$$
|1+\EIGHT\eta^2\gamma ^2-\HALF\gamma |\,\le\,\gamma ~.
$$
We
polarize the {\em second} integral in \equ{split} (but not the first)
and bound it (in modulus) by
$$
\eqalign{|1+\EIGHT\eta^2\gamma ^2-\HALF\gamma |&\idx  \hdd (|\gu_t \gv_t |+|\gu_t '\gv_t '|)\,\le\,
\gamma \idx  \hdd (|\gu_t \gv_t |+|\gu_t '\gv_t '|)\cr
\,&\le\,
\idx  \hdd \bigg(8\gamma  ^2\bigl(\gu_t ^2+(\gu_t ')^2\bigr)+\EIGHT
\bigl(\gv_t ^2+(\gv_t ')^2\bigr)\bigg) ~.}
\EQ{uvbound}
$$
Combining \equ{hdd2} with the decomposition \equ{split} and the bound
\equ{uvbound}, we find 
$$
\eqalign{
\HALF\partial _t J(\gu_t,\gv_t) \,\le\,&-\idx  \hdd \bigl(\gv_t ^2+(\gv_t ')^2\bigr
)\bigl(1-\HALF\eta^2\gamma -\CX{dd}\delta-\textstyle{1\over 8}\bigr)\cr
&-\HALF\gamma \idx  \hdd
\bigl(\gu_t ^2+2(\gu_t ')^2+(\gu_t '')^2\bigr)\bigl(1-{\CX{dd}\delta\over
\gamma }-2\CX{dd}\delta \bigr)\cr
&+\idx  \hdd \bigl(\gu_t ^2+(\gu_t ')^2\bigr)
\bigl((8\gamma ^2+\HALF\CX{dd}\delta\gamma )+\HALF\gamma (1-{\CX{dd}\delta\over
\gamma }-2\CX{dd}\delta )\bigr )\cr 
&-\EIGHT\eta^2\gamma ^2\idx  \hdd\bigl(\gu_t \gv_t +\gu_t '\gv_t '\bigr)~.\cr
}
\EQ{hdd3}
$$
The bizarre decomposition of the terms involving $(\gu_t ')^2$ will become
clear soon.
Note that by our choice of constants, \equ{hdd3} can be simplified to
the slightly less good bound
$$
\eqalign{
\HALF\partial _t J(\gu_t,\gv_t) \,\le\,&-\HALF\idx  \hdd \bigl(\gv_t ^2+(\gv_t ')^2\bigr
)\cr
&-\FOUR\gamma \idx  \hdd
\bigl(\gu_t ^2+2(\gu_t ')^2+(\gu_t '')^2\bigr)\cr
&+16\gamma ^2\idx  \hdd \bigl(\gu_t ^2+(\gu_t ')^2\bigr)
\cr 
&-\EIGHT\eta^2\gamma ^2\idx  \hdd\bigl(\gu_t \gv_t +\gu_t '\gv_t '\bigr)~.\cr
}
\EQ{hdd3a}
$$
We project onto high momenta, and exploit the contraction
properties:
{\em We assume from now on that $\QQ_{k_*}\gu_0 =0$ and $\QQ_{k_*}\gv_0 =0$.}
Note that if this property holds at time zero, it holds for all times
for the evolution defined by $\LL$, because $\LL$ commutes with
$\QQ_{k_*}\oplus \QQ_{k_*}$.
Using the bounds of \sec{high}, we have
$$
\eqalign{
\idx  \hdd \gu_t ^2 \,&\le\, \nu \idx  \hdd (\gu_t ')^2~,\cr
\idx  \hdd (\gu_t ')^2 \,&\le\, \nu \idx  \hdd (\gu_t '')^2~,\cr
}\EQ{nubounds}
$$
where
$$
\nu\,=\,\CX{nu}k_*^{-2}~.
$$
Thus, \equ{hdd3a} can be improved to
$$
\eqalign{
\HALF\partial _t J \,\le\,&
-\HALF\idx  \hdd \bigl(\gv_t ^2+(\gv_t ')^2\bigr
)\cr
&-\FOUR\gamma \idx  \hdd
\bigl(\gu_t ^2+2(\gu_t ')^2+(\gu_t '')^2\bigr)\cr
&+16\nu\gamma ^2\idx  \hdd\bigl( (\gu_t ')^2+(\gu_t '')^2\bigr )\cr
&-\EIGHT\eta^2\gamma ^2\idx  \hdd\bigl(\gu_t \gv_t +\gu_t '\gv_t '\bigr)~.\cr}
\EQ{hdd4}
$$
This leads to a bound of the form
$$
\eqalign{
\HALF \partial _t J\,\le\,& - \HALF\eta^{-2} \idx  \hdd
\eta^2\bigl(\gv_t ^2+(\gv_t ')^2\bigr)\cr
&-\EIGHT\gamma \idx  \hdd\bigl(\gu_t ^2+2(\gu_t ')^2+(\gu_t '')^2\bigr) \cr
&-\EIGHT\eta^2\gamma ^2\idx 
\hdd\bigl(\gu_t \gv_t +\gu_t '\gv_t '\bigr)~.\cr
}
$$
Taking the least good bounds above, we finally get the decay of the
high frequency part (since $\eta^{-2}\ge\gamma $):
$$
\eqalign{
\HALF \partial _t J\,&\equiv\,\HALF\partial_t\idx  \hdd \bigg(
\eta^2\bigl(\gv_t ^2+(\gv_t ')^2\bigr)+\bigl(\gu_t ^2+2(\gu_t ')^2+(\gu_t '')^2\bigr )+\eta^2\gamma 
\bigl(\gu_t \gv_t +\gu_t '\gv_t '\bigr)\bigg)\cr
\,&\le\, - \EIGHT\gamma  \idx  \hdd \bigg(
\eta^2\bigl(\gv_t ^2+(\gv_t ')^2\bigr)+\bigl(\gu_t ^2+2(\gu_t ')^2+(\gu_t '')^2\bigr )+\eta^2\gamma 
\bigl(\gu_t \gv_t +\gu_t '\gv_t '\bigr)\bigg)\cr
\,&=\,-\EIGHT\gamma  J\,=\,-{\textstyle{1\over 2560}} \min(\eta^{-2},k_*^2/\CX{nu})J~.\cr
}
\EQ{enfin}
$$

Thus we have shown the 
\CLAIM{Proposition}{J}There is a (small) $\eta_0>0$ such that for all
$\eta<\eta_0$ the following holds for
the functional 
$$
J_{\hdd,\gamma }(\gu_t,\gv_t)\,=\,
\idx  \hdd(x)\bigg(\eta^2\bigl(\gv_t ^2+(\gv_t ')^2\bigr)
+\gu_t ^2+2(\gu_t ')^2+(\gu_t '')^2+\eta^2\gamma \bigl(\gu_t \gv_t +\gu_t '\gv_t ')\bigr)\bigg)(x)~.
$$
Let $(\gu_t,\gv_t)=e^{\LL t}(\gu_0,\gv_0)$, and assume
$(\gu_0,\gv_0)\in\KK_{k_*}\oplus\KK_{k_*}$.
Then
$$
J_{\hdd,\gamma }(\gu_t,\gv_t)\,\le\, \exp\bigl(-\gamma t/80\bigr) \cdot
J_{\hdd,\gamma }(\gu_0,\gv_0)~,
\EQ{jdef2}
$$
where
$$
\gamma \,=\,\min(\eta^{-2},k_*^2/\CX{nu})/320~.
\EQ{gammadef} 
$$

We come now back to the definition \equ{jdef2} of $J$, and compare it
to the norm $\|\cdot\|_{\hdd,2}$ defined in Eq.\equ{ldd}.
These two quantities define equivalent topologies when considered on
$\KK_{k_*}\oplus \KK_{k_*}$.
\CLAIM{Lemma}{uv}On $\KK_{k_*}\oplus \KK_{k_*}$ one has the bound
$$
\eta^2\gamma \left |\idx  \hdd\bigl(\gu \gv +\gu '\gv '\bigr) \right| 
\,\le\,\idx  \hdd \bigg(\HALF\bigl((\gu ')^2+(\gu '')^2\bigr)+\EIGHT 
\bigl(\gv ^2+(\gv ')^2\bigr)\bigg)
~.
$$

\REMARK This lemma eliminates the somewhat arbitrary quantity $\gamma $ from the
topology, see \clm{JL} below.
\PROOF This is a combination of earlier estimates.
Indeed, we have already seen in Eq.\equ{uvbound} that
the mixed terms in Eq.\equ{split} can be bounded by
$$
\eta^2\gamma \left |\idx  \hdd\bigl(\gu \gv +\gu '\gv '\bigr) \right| 
\,\le\,\eta^2\idx  \hdd \bigg(8\gamma  ^2\bigl(\gu ^2+(\gu ')^2\bigr)+\EIGHT
\bigl(\gv ^2+(\gv ')^2\bigr)\bigg)\,\equiv\,X~.
$$
Furthermore, by \equ{nubounds} and the choice of $k_*$, we get
$$
X\,\le\,\eta^2\CX{nu}k_*^{-2}\idx  \hdd \bigg(8\gamma 
^2\bigl((\gu ')^2+(\gu '')^2\bigr)+\EIGHT 
\bigl(\gv ^2+(\gv ')^2\bigr)\bigg)~.
$$
Since we have also chosen $\gamma \,=\,\min(\eta^{-2},k_*^2/\CX{nu})/320$,
we get
finally
$$
X\,\le\,\idx  \hdd \bigg(\HALF\bigl((\gu ')^2+(\gu '')^2\bigr)+\EIGHT 
\bigl(\gv ^2+(\gv ')^2\bigr)\bigg)
~,
$$
which is what we asserted.\QED

Recall the definition \equ{pdef} of
the
projection $\PP_{k_*}$ onto momenta larger than $k_*$.
From \clm{uv} and \clm{J} we have immediately, with the notation of
\equ{ldd} and \equ{jdef0} in the topology of $\HH_{\delta,\loc}^2$ (which
does not depend on $\delta_0$).

\CLAIM{Theorem}{JL}Assume $\eta_0$ and $k_*$ satisfy \equ{eta0} and
\equ{defk*}, and assume  $\delta\le1/(40\CX{dd})$.
For all $\eta$
satisfying 
$0<\eta<\eta_0$ the following holds:
If $\|(\gu _0,\gv _0)\|_{\hdd,2}<\infty $ 
and  $(\gu _t,\gv _t)=e^{\LL t}(\gu _0,\gv _0)$
then one has the bounds
$$
J_{\hdd,\gamma }(\PP_{k_*}\gu _t,\PP_{k_*}\gv _t)/2 \,\le\,\|(\PP_{k_*}\gu _t,\PP_{k_*}\gv _t)\|_{\hdd,2}^2\,\le\,2J_{\hdd,\gamma }(\PP_{k_*}\gu _t,\PP_{k_*}\gv _t)~,
\EQ{JL}
$$
and
$$
\|(\PP_{k_*}\gu _t,\PP_{k_*}\gv _t)\|_{\hdd,2}\,\le\,4\exp(-\gamma t/80)\cdot\|(\PP_{k_*}\gu _0,\PP_{k_*}\gv _0)\|_{\hdd,2}~,
\EQ{Lt}
$$
where
$\gamma =\min(\eta^{-2},k_*^2\CX{nu})/320$.

\SECT{full}{The evolution of differences}

In this section, we combine the results of Sections
\lab{s}{linearized}  and \lab{s}{highfreq} into bounds on
the evolution of the difference of solutions to \equ{problem}. 
We will first treat the general case, and show a bound which diverges
exponentially with time, and then we will treat the high frequency
case where we have decay.
We consider two initial conditions, and their respective evolutions
under the semi-flow defined by \equ{problem}. We call these functions
$(u_1,v_1)$ and $(u_2,v_2)$, respectively.
The evolution for the difference
$(\gu ,\gv )=(u_1-u_2,v_1-v_2)$
takes now the form 
$$
\eqalign{
\dot \gu \,&=\,\gv ~,\cr
\eta^2\dot \gv \,&=\,-\gv +\gu ''+M(u_1,u_2)\gu ~,
\cr
}
\EQ{diff1a}
$$
where $M(u_1,u_2)$ is defined by
$M(u_1,u_2)(u_1-u_2)=U'(u_1)-U'(u_2)$.
It will be
useful to rewrite this system of equations as
$$
\V{\dot \gu_t}{\dot \gv_t}\,=\,\LL\V{\gu_t}{\gv_t}+
\V{0}{M_{u_{1,t},u_{2,t}}\gu_t}~.
\EQ{llforma}
$$
Note that $M(u_1,u_2)$ is really a space-time dependent coefficient of the
linear problem \equ{llform}.
The important observation is now that on the attracting set, {\it
i.e.}, for all sufficiently large $t$ (depending on the initial
conditions $u_1, u_2, v_1, v_2$) we have,
by \clm{f1}, a universal bound 
$$
\sup_{x\in \real}|M(u_{1,t}(x),u_{2,t}(x))|+\sup_{x\in \real}|\partial_xM(u_{1,t}(x),u_{2,t}(x))|\,\le\,M_*~.
\EQ{mmm}
$$

Since we already know bounds on the solution, we can write it as
follows for $\gu_t(x)=u(x,t)$ and $\gv_t(x)=v(x,t)$:
$$
\V{\gu_t}{\gv_t}\,=\,e^{t\LL} \V{\gu_0 }{\gv_0 }+\int _0^t\d{s} e^{(t-s)\LL}
\V{0}{M_{u_{1,s},u_{2,s}}\gu_s}~.
\EQ{repr1} 
$$

\CLAIM{Proposition}{expo1}Assume $(u_{1,0},v_{1,0})$ and
$(u_{2,0},v_{2,0})$ are in $\GG$. Let $\gu_t=u_{1,t}-u_{2,t}$
and let  $\gv_t=v_{1,t}-v_{2,t}$. There are constants $\C{MM1}$ and
$\C{MM2}$
such
that for all $t>0$,
$$
\|(\gu_t,\gv_t)\|_{\hdd,2}\,\le\,
\CX{MM1}e^{\CX{MM2}t}\|(\gu_0,\gv_0)\|_{\hdd,2}  
~.
\EQ{MM1} 
$$

\PROOF We have already seen in \equ{mmm} that $|M(u_{1,t},u_{2,t})|$
and its derivative
are bounded and then the result follows at once from the representation 
\equ{repr1}  and the bound of \clm{frakuv}.\QED

The handling of the high frequency part $\PP_{k_*}(\gu_t,\gv_t)$ is
similar.
Instead of \equ{llforma}, we get
$$
\V{\partial _t\PP_{k_*}\gu_t}{\partial
_t\PP_{k_*}\gv_t}\,=\,\LL\V{\PP_{k_*}\gu_t}{\PP_{k_*}\gv_t}+
\V{0}{\PP_{k_*}M_{u_{1,t},u_{2,t}}\gu_t}~.
\EQ{llformb}
$$ 
The solution of this problem is
$$
\V{\PP_{k_*}\gu_t}{\PP_{k_*}\gv_t}\,=\,e^{t\LL} \V{\PP_{k_*}\gu_0
}{\PP_{k_*}\gv_0 }+\int _0^t\d{s} e^{(t-s)\LL}
\V{0}{\PP_{k_*}M_{u_{1,s},u_{2,s}}\gu_s}~.
\EQ{look} 
$$
What is important here is that in {\em both} terms the operator $\LL$
acts on functions with {\em high} frequencies.

\CLAIM{Proposition}{expo2}Assume $(u_{1,0},v_{1,0})$ and
$(u_{2,0},v_{2,0})$ are in $\GG$. Let $\gu_t=u_{1,t}-u_{2,t}$
and let  $\gv_t=v_{1,t}-v_{2,t}$. There are constants $\C{MM9}$, $\C{MM3}$, and
$\C{MM4}$
such
that for all $t>0$,
$$
\|(\PP_{k_*}\gu_t,\PP_{k_*}\gv_t)\|_{\hdd,2}\,\le\,
\bigg(\CX{MM9}e^{-\gamma t/80}+\CX{MM3}{e^{\CX{MM4}t}\over
\gamma }\bigg)\|(\gu_0,\gv_0)\|_{\hdd,2}   
~,
\EQ{MM2} 
$$
where $\gamma =\min(\eta^{-2},k_*^2/\CX{nu})/320$.

\REMARK In fact, one can choose $\CX{MM4}=\CX{MM2}$.
\PROOF We use again \equ{mmm} to bound $M$ and $\partial _xM$. 
Furthermore, $\PP_{k_*}$ is bounded and then the result
follows at once from the representation  
\equ{look} and the bound \equ{Lt} of \clm{JL} for the first term of
\equ{MM2} and 
additionally 
the bound \equ{MM1} of \clm{expo1} for the second.\QED

\SECT{bounds}{Covering the attracting set}

We define a new norm by
$$
\|(u,v)\|_{\delta,L,2}\,=\,\sup_{\xi\in[-L,L]}
\|(u,v)\|_{h_{\delta,\xi},2}~,
\EQ{newnorm}
$$
where
$$
h_{\delta,\xi}(x)\,=\,{1\over (1+\delta^2(x-\xi)^2)^2}~,
$$
and 
$\|(u,v)\|_{\hdd,2}$ was defined in \equ{ldd}. 
This norm looks at a ``window'' of size $2L$ in 
$\HH_{\delta,\loc}^2$.
For $\epsilon >0$ we define $N_L(\epsilon )$ as the minimum number of
balls of radius $\epsilon $ (in the norm $\|\cdot\|_{\delta,L,2}$),
needed to cover the attracting set $\GG$.
\CLAIM{Theorem}{balls}There exist finite constants $A$, and $\C{ball}$
such that for all $\epsilon$ with $0<\epsilon <1$ and all
$L>A/\epsilon $ one has
$$
N_{L-A/\epsilon }(\epsilon )\,\le\,\CX{ball}^L N_L(2\epsilon
)~.
\EQ{balls}
$$ 

\PROOF We denote $t\mapsto\Phi^t$ the flow defined by \equ{problem}.
Let $\BB$ be
a finite collection of balls
of radius $\epsilon $ in the topology defined by
$\|\cdot\|_{\delta,L,2}$ which cover $\GG$.

We next define a natural {\em unit of time}, $\tau _*$.
We recall the definition \equ{gammadef} of $\gamma$: 
$\gamma =\min(\eta^{-2},k_*^2/\CX{nu})/320$.
We define
$$
\tau_*\,=\,{b\over \gamma}\log \gamma~,
\EQ{taudef}
$$
where the (small positive) constant $b$ is chosen such that the factor
in \equ{MM2} is minimal and when $\gamma$ is large (say,
$\gamma>\gamma_0$), we get 
$$
\CX{MM9}e^{-\gamma \tau _*/80}+\CX{MM3}{e^{\CX{MM4}\tau _*}\over
\gamma }\,\le\, \gamma^{-\kappa}~,
\EQ{theineq}
$$
for some $\kappa>0$. We will use this bound in the sequel.

Since the flow 
$\Phi^t$ leaves $\GG$ invariant,
we see that 
$$
\Phi^{t+\tau}(\GG)\,\subset\,\bigcup_{B\in\BB} \Phi^\tau \bigl(B \cap
\Phi^t(\GG)\bigr )\,=\,\bigcup_{B\in\BB} \Phi^\tau \bigl(B \bigr )~.
$$
Consider now any of the $B$ in $\BB$. We are going to cover
$\Phi^{\tau_*}(B)$ by balls of radius $\epsilon /2$. Let $\fff_0$ and $\ggg_0$
be
two elements of the ball $B$ and assume furthermore $\fff_0$ and $\ggg_0$
are in the attracting set $\GG$. This implies
$$
\|\fff_0-\ggg_0\|_{\delta ,L,2}\,\le\,\epsilon ~,
\EQ{bound1}
$$
and, since $\fff_0$ and $\ggg_0$ are in the global attracting set $\GG$,
we also
have
$$
\|\fff_0-\ggg_0\|_{\hdd,2}\,\le\,\C{initial}~,
\EQ{bound2}
$$
for some constant $\CX{initial}$.
With ${\tau_*}$ as in \equ{taudef}, we let
$$
\fff \,=\,\Phi^{\tau_*} (\fff_0)~,\quad \ggg\,=\,\Phi^{\tau_*} (\ggg_0)~.
$$
We then rewrite  $\fff-\ggg$ as
$$
\fff-\ggg\,=\,\PP_{k_*}(\fff-\ggg) +\QQ_{k_*}(\fff-\ggg)~,
\EQ{decomp1}
$$
where (the direct sums of) $\PP_{k_*}$ and $\QQ_{k_*}$ are the high- and low-momentum
projections introduced earlier (in \equ{Qdef2} and \equ{pdef}). 
Our aim is to bound this difference in the norm
$\|\cdot\|_{\delta ,L-A/\epsilon,2 }$, where $A$ is a large
constant to be determined later. 

We begin with $\PP_{k_*}(\fff-\ggg)$.
By our choice of $\tau_*$ in \equ{taudef}, we have, by \equ{theineq}
and
\clm{expo2}, 
$$
\big\|\PP_{k_*}\bigl(\fff-\ggg\bigr
)\big\|_{\delta,L-A/\epsilon ,2 }\,\le\,\big\|\PP_{k_*}\bigl(\Phi
^{\tau_*} (\fff_0)-\Phi ^{\tau_*} (\ggg_0)\bigr 
)\big\|_{\delta,L ,2 }\,\le\,
\gamma^{-\kappa} \epsilon ~.
$$
{\em We now fix $\eta_0>0$ so small and $k_0$ so large
(and at least satisfying Eqs.\equ{eta0} and \equ{defk*})
such that for all $\eta<\eta_0$ and all $k_*>k_0$ one has}
$$
\gamma^{-\kappa}\,=\,\bigl(\min(\eta^{-2},k_*^2/ \CX{nu})/320\bigr
)^{-\kappa}\,\le\, {1\over 8}~,
\EQ{finaleta0}
$$
and also $\gamma>\gamma_0$, see \equ{theineq}.
Summarizing the bounds for this piece, we get
$$
\|\PP_{k_*}(\fff-\ggg)\|_{\delta ,L-A/\epsilon ,2}\,\le\,{\epsilon \over
8}~.
\EQ{final1}
$$

We bound next $\QQ_{k_*}(\fff-\ggg)$ by decomposing it into a finite sum plus
a remainder. 
We will work with the two components of $\QQ_{k_*}\fff$  or
$\QQ_{k_*}\ggg$ separately. Since the norm on the first component has
2 derivatives and the norm on the second only 1, we will deal only
with the first case and leave the other case to the reader.

We will work with the 
notion of Bernstein class $\BB_{\sigma}(K)$, defined by
$$
\BB_{\sigma}(K)\,=\,\bigg\{ h~\bigg|~|h(x+iy)|\le
Ke^{\sigma|y|}{\text{ for all }} x,y\in\real\biggr\}~.
$$ 
If $h\in\BB_{\sigma}(K)$, it can be represented by the Cartwright
interpolation formula [KT, Eq. (191)] (or [B] for a proof)
with $\sigma'=\pi/2$ and
$\omega=\pi/4$) as
$$
h(x)\,=\,{\sin(2 \sigma x)\over 4} \sum_{j=-\infty }^\infty 
(-1)^j {\sin(\cartarg)\over (\cartarg)^2}h(x_j)~,
\EQ{cart}
$$
where the $x_j= {j \pi\over 2\sigma} $ are discrete sampling
points. This class is useful in our context because of
\CLAIM{Lemma}{zzz}There is a constant $\C{BB}$ such that if
$u\in\L^\infty $, then 
$$
\QQ_{k_*}u\in \BB_{2k_*}(\CX{BB}k_*\|u\|_\infty) ~.
\EQ{incomplete} 
$$

\PROOF This amounts to saying that a function with frequency support
in 
$[-2k_*,2k_*]$ is in the Bernstein class. This is almost obvious, except for
the smooth cutoff.
In fact, with the function $\theta$ as defined in \sec{high}, we consider
$$
\int \d{k} e^{ik(x+iy)}\theta(k/k_*)\,=\,
k_*\int \d{\ell} e^{ik_*\ell (x+iy)}\theta(\ell)~,
\EQ{exponential} 
$$
which is in $\L^1(\,\,\d{x})$ for any $y\in\real$. And the
$\L^1(\,\,\d{x})$ norm is bounded by
$$
\OO(1)\bigl(k_*+k_*^{-1}\bigr)e^{2k_* |y|}~.
$$
Therefore, the convolution operator defined by \equ{exponential} maps
$u$ to $\BB_{2k_*}(\OO(k_*)\|u\|_\infty )$.\QED

We next bound the functions appearing in
\equ{cart} in our favorite topology:
\CLAIM{Lemma}{cart}Let $\sigma>2$ and let $f_j$ be defined by
$$
f_j(x)\,\equiv\,{\sin\bigl (4(\cartarg)\bigr )\sin(\cartarg)\over
4(\cartarg)^2}\,=\,(-1)^j{\sin\bigl (2\sigma x\bigr )\sin(\cartarg)\over
4(\cartarg)^2}~. 
$$
There is a constant $\C{bound}$ independent of $j$ and $\xi$, 
such that for all
$j$ and $\xi$
one has:
$$
\int\dx h_{\delta}(x-\xi)\bigg( f_j^2(x)+2\bigl(f_j'(x)\bigr
)^2+\bigl(f_j''(x)\bigr )^2\bigg) \,\le\,{\sigma^4\CX{bound}\over
1+(2\sigma\xi-\pi j)^4}~. 
\EQ{thebound}
$$

\REMARK The numerical coefficient $\CX{bound}$ depends on $\delta$,
but $\delta $ has been fixed in Eq.\equ{deltadef}:
$\delta=1/(40\CX{dd})$.
\PROOF The function $f_j$ can be bounded as
$$
|f_j(x)|\,\le\, \left|{\sin\bigl (4(\cartarg)
\bigr )\sin(\cartarg)\over 4(\cartarg)^2}\right
|\,\le\, {\C{x4}\over 1+ (\cartarg)^2 }~,
$$
since the numerator vanishes simultaneously with the denominator (and
to order 2). Similarly, the derivative is bounded by
$$
\left|\partial_x^\ell  \left({\sin\bigl (4(\cartarg\bigr )
\sin(\cartarg)\over 4(\cartarg)^2}\right
)\right|\,\le\,{\C{x5}(\sigma/2)^\ell\over  1+ 4(\cartarg)^2}~,\quad\ell=1,2~,
\EQ{der}
$$
since $\sigma>2$ by assumption.
It follows that
$$
\int\dx h_{\delta,\xi}(x) |f_j(x)|^2\,\le\,
\C{x44}\int\dx {1\over (1+\delta^2(x-\xi)^2)^2}\,\cdot\,
{1\over (1+(\cartargfour)^2)^2}~.
$$
Setting $\rho=\min(\delta,2\sigma)$, we find that this is bounded by
$$
\int\dx h_{\delta,\xi}(x) |f_j(x)|^2\,\le\,{\C{x99}\over \rho} {1\over
\bigl (1+\rho^2 (\xi-{\pi j\over 2\sigma})^2\bigr )^2}~.
$$
In view of \equ{der} one gets a similar bound for the derivatives, and
thus \equ{thebound} follows.\QED

Consider the element $(u,v)\in\GG$. We know that
$\|(u,v)\|_{\delta,\loc,2}\le\CX{77}$. 
For the first component, $u$, this means
$$
\sup_{\xi\in\real}\int\dx h_\delta(x-\xi) \bigg( |u(x)|^2 +
2|u'(x)|^2+ |u''(x)|^2\bigg)\,\le\,\CX{77}^2~.
$$
From this, we conclude using the Sobolev inequality in the form of
\clm{Sobolev} that
$\|u\|_\infty \le\C{h22}$ for 
some finite $\CX{h22}$. By \clm{zzz} we then get that
$\|\QQ_{k_*}u\|_\infty \le \CX{BB}k_*\CX{h22}$ and furthermore,
$\QQ_{k_*}u\in\BB_{2k_*}(\CX{BB}k_*\CX{h22})$. Thus, we can apply the
Cartwright formula to $h=\QQ_{k_*}u$, with $\sigma=2k_*$.

Throughout,
$L k_*$ has to be sufficiently large.
We define
$$
\eqalign{
\SS_{L}(h)\,&=\,{\sin(4 \kstar x)\over 4} \sum_{|j|\le
2L\kstar}  
(-1)^j {\sin(\cartk)\over (\cartk)^2}h(x_j)~,\cr
}\EQ{S}
$$
where $x_j={j\pi\over 4\kstar}$ are the discrete sampling points.
We decompose
$$
\QQ_{k_*}u\,=\,\Big(h-\SS_{L}(h)\Big)+\SS_{L}\bigl(h\bigr )~.
\EQ{sh1} 
$$
The first term in \equ{sh1} will be small because 
$h-\SS_{L}(h)$ is the remainder
of the converging sum in \equ{cart}, and for the second one we will
use a covering argument.

We first show that $X_L\equiv h-\SS_{L}(h)$ is small when $L$ is large.
The difference can be written as
$$
\bigl(h-\SS_L (h)\bigr)(x)\,=\,
\sum_{|j|>2L\kstar}(-1)^j {\sin(4 \kstar
x)\sin(\cartk)\over 4(\cartk)^2}h({j\pi\over 4k_*}) ~.
$$
Using \equ{thebound}, we get as a bound for $X_L$ {\em when $|\xi|\le L$}:
$$
\eqalign{
\bigg(\int \dx \hdd(x-\xi)&\bigl(|X_L(x)|^2+2|X_L'(x)|^2+|X_L''(x)|^2\bigr)
\bigg)^{1/2}\cr\,&\le\,\CX{bound} \sum_{|j|\ge
2L\kstar}{1\over  1+(k_*\xi -{\pi j\over 4})^2}\,\le\,{\C{n1}\over
1+\bigl|L-|\xi|\bigr|}~.  }
$$ 
This argument can be repeated for the second component.
Since in the definition \equ{newnorm} of $\|\cdot\|_{\delta
,L-A/\epsilon ,2}$  we have $|\xi|\le L-A/\epsilon $,
we find a bound on the exterior part of $\QQ_{k_*}(\phi-\psi)$:
$$
\Big\|\QQ_{k_*}(\phi-\psi)-\SS_{L }\bigl(\QQ_{k_*}(\phi-\psi)\bigr
)\Big\|_{\delta ,L-A/\epsilon ,2}\,\le\,{\C{L}\epsilon\over A
}\sup_{\xi\in\real}\|\phi-\psi\|_{h_{\delta ,\xi},2}\,\le\,{\CX{L}\epsilon\over A
}\CX{initial}~,
\EQ{for9}
$$
using \equ{bound2}. Clearly, if $A$ is sufficiently large (but
independent of $\epsilon $ and $L$), we get the bound
$$
\Big\|\QQ_{k_*}(\phi-\psi)-\SS_{L }\bigl(\QQ_{k_*}(\phi-\psi)\bigr
)\Big\|_{\delta ,L-A/\epsilon ,2}\,\le\,{\epsilon\over 8}~.
\EQ{part1}
$$

We finally deal with the central part, namely
$\SS_{L}(\QQ_{k_*}(\phi-\psi))$.
This is described in
\CLAIM{Proposition}{finite}There is a constant
$\C{L2}$ such that the 
following holds.
Let $B$ be a ball of radius $\epsilon  $ in the topology defined by
$\|\cdot\|_{\delta ,L,2}$.
Then the set $\SS_{L}(B\cap\GG)$ can be covered by no more than
$$
\CX{L2}^L 
$$
balls of radius $\epsilon/8$. 

\PROOF Since $\phi,\psi\in \GG$, \clm{zzz} implies
$\QQ_{k_*} (\phi-\psi)\in \BB_{{2k_*}}(X)$,
where 
$$
X\,=\,{\rm diam}_{\L^\infty } (\GG) \CX{BB}k_*~.
$$
Moreover, from \clm{Qa} we deduce
$$
\|\QQ_{k_*} (\phi-\psi)\|_{\delta,L,2}\,\le\,\C{p1}\epsilon ~.
$$
Using the Sobolev inequality from \clm{Sobolev}, this implies
$$
\sup_{x\in[-L,L]}\bigl|\bigl(\QQ_{k_*}(\phi-\psi)\bigr
)(x)\bigr|\,\le\,\C{p2}\epsilon ~.
$$
We use next the bounds
$$ 
\bigl|\bigl(\QQ_{k_*}(\phi-\psi) \bigr )(x_{j})\bigr |\,\le\,\CX{p2} \epsilon
~,
$$
for $|j|< {\lpoint}$. We let $n$ be a large integer which
will be fixed at the end of the proof.
The set of values of each of the 2 components of
$\bigl(\QQ_{k_*}(\phi-\psi) \bigr 
)(x_{j})$ can be covered by $8n\CX{p2}$ balls of radius
$\epsilon /(4n)$, for each of the $2(2L\kstar)+1$ possible values of
$j$. We bound now in detail the sum in
$\SS_{L}\bigl(\QQ_{k_*}(\phi-\psi) \bigr )$ as defined in \equ{S}.

We fix one of the  $(8n\CX{p2})^{4(\lpoint)+1}$ grid points
for the components of $\bigl(\QQ_{k_*}(\phi-\psi) \bigr
)(x_{j})$. For each component, we get a set of $2(\lpoint)+1$ numbers
$q_\ell$, with $|\ell|\le \lpoint$.
We pick numbers $r_\ell$ satisfying $|r_\ell-q_\ell|<\epsilon
/(4n)$ for all 
$\ell$ and we want to show that the function
$$
\Delta(x)\,=\,{\sin(4\kstar x)\over 4} \sum_{|j|\le
2L\kstar}  
(-1)^j {\sin(\cartk)\over (\cartk)^2} (r_j-q_j) 
$$
has a $\|\cdot\|_{\delta,L,2}$ norm less than $\epsilon /8$. This will
clearly suffice to show \clm{finite}.
By \clm{cart}, we get
$$
\|\Delta\|_{\delta ,L,2}\,\le\,\C{delta}\sup_{|\xi|\le L}
\sum_{|j|\le
2L\kstar}  {1\over 1+(4\kstar\xi-\pi j)^2} {\epsilon \over
4n}\,\le\,\C{delta2}{\epsilon \over n}~.
$$
We choose $n=8\CX{delta2}$,
and we see that,  all in all, one needs
$(8n\CX{p2})^{4(\lpoint)+1}\le\CX{L2}^L$ balls of radius $\epsilon
/8$ to cover $\SS_{L}(B\cap\GG)$. (Note that $\CX{delta}$ and $\CX{delta2}$
depend on $k_*$. In fact they are bounded by $\OO(k_*^4)$.)\QED

\LIKEREMARK{Proof of \clm{balls}}We combine now the various estimates
to prove \equ{balls}. Let $B$ be one of the $N_L(\epsilon )$ balls of
radius $\epsilon $ needed to cover $\GG$ and let $f\in B\cap\GG$.
All we need to show is that the set of all $g\in B\cap\GG$ can be
covered by $\CX{ball}^L$ balls of radius $\epsilon/2$ in the topology
of the norm $\|\cdot\|_{\delta ,L-A/\epsilon ,2}$. We decompose $\phi-\psi$
according to \equ{decomp1} and then $\QQ_{k_*}(\phi-\psi)$ according to
\equ{sh1}, so that we have three terms.
The first is bounded by $\epsilon /8$ using \equ{final1} and the
second is bounded by \equ{part1}. Since $\SS_{L}(B\cap\GG)$ can be
covered by $\CX{L2}^L$ balls of radius $\epsilon/8 $ in the norm
$\|\cdot\|_{\delta ,L,2}$ it can also be covered by the same number of
balls in the norm $\|\cdot\|_{\delta ,L-A/\epsilon ,2}$. Thus the sum
of the three contributions can be covered by $\CX{L2}^L$ balls of
radius $3\epsilon /8<\epsilon /2$. The proof of \clm{balls} is complete.\QED

\SECT{entropy}{The $\twelvepoint\epsilon $-entropy and the topological entropy}
\SUBSECTION{Introduction}

\rm In this section, we exploit the results obtained so far to show that
the $\epsilon$-entropy and the 
topological entropy per unit length can be defined for the
Eq.\equ{problem}. The reasoning here is very close to the one used in
[CE2], and so there is no need to repeat it here. What needs however
some special attention is the choice of topology for which the entropy
per unit length can be defined. We basically need a topology which has
a {\em submultiplicativity property} which we define below. The most
simple example of such a topology was used in [CE2], namely $\L^\infty
$.
The property which we used there is that if a set $S$ of functions is
defined on the union of 2 adjacent intervals, say $I_1\cup I_2$, then
the following is true: If $S$ restricted to $I_1$ can be covered by
$N_{I_1}$ balls of radius $\epsilon $ in $\L^\infty(I_1) $, and $S|_{I_2}$ can
be covered by
$N_{I_2}$ balls in $\L^\infty (I_2)$, then $S|_{I_1\cup I_2}$ can be covered by $N_{I_1}\cdot
N_{I_2}$ balls in $\L^\infty (I_1\cup I_2)$ (all of radius $\epsilon $).
In $\L^\infty $, this property is obvious:
Let $B_{1,i}$, with $i=1,\dots,N_{I_1}$ be the balls covering $S|_{I_1}$
and $B_{2,j}$, with $j=1,\dots,N_{I_2}$ those covering $S|_{I_2}$. Then
one can just take the set $S_{i,j}$ of functions
$$
S_{i,j}\,=\,\bigg\{ f~\bigg | ~ f|_{I_1}\in B_{1,i} \text{ and }
 f|_{I_2}\in B_{j,2}\bigg\}~,
$$
and this {\em is} a ball of radius $\epsilon $ in $\L^\infty (I_1\cup
I_2)$.

The difficulty with topologies which are finer than $\L^\infty $ is
that we have to patch the functions on $I_1$ and $I_2$ together in
such a way that the patched function is an element of a ball in the
topology on $I_1\cup I_2$. We do not know how to do this in the
topologies used in the earlier sections, and therefore we go to a new
topology in which the submultiplicativity property holds in the sense that
there is a constant $C=C(\epsilon )$ {\em independent of $I_1$ and
$I_2 $} such that the 
functions on the union of $I_1$ and $I_2$ can be covered by 
$$
N_{I_1}(\epsilon )\cdot N_{I_2}(\epsilon )\cdot C (\epsilon )
\EQ{nnn}
$$
balls of radius $\epsilon $. It is well known from the literature on
statistical mechanics (see {\it e.g.}, Ruelle [R]) and easy to see that
this weaker form of submultiplicativity suffices to prove the existence of
limits (of the logarithms) per unit length.

The topology we will use is $\WWW$, defined by
$$
\| f\|_{\WWW}\,\equiv\,
\max\bigl(\sup_{x\in \real }|f(x)|,~ \sup_{y\in \real} |f'(y)| \bigr
)~.
\EQ{wwwdef}
$$
This is a ``good'' topology for our problem, because we can control the
evolution of functions in $\WWW$. However, it is obvious that the
submultiplicativity property is not immediate, since the matching of
functions has to be continuous and once differentiable.

\SUBSECT{subadditivity}{Submultiplicativity in $\WWW$}

We develop here the estimates leading to Eq.\equ{nnn} for balls in
$\WWW$.
Our main result will be \clm{both}.
We let $\R>5$ be a large constant which will be determined in Eq.\equ{r}
below.
\LIKEREMARK{Notation}Throughout, we will use the notation
$$
|g|_I \,=\,\sup_{x\in I} |g(x)|~.
$$
We let $\WWW_I$ be the space of continuously
differentiable functions $g:I\to\real$, equipped with the norm
$$
\|g\|_I\,=\,\max \bigl(|g|_I,|g'|_I\bigr)~.
$$
(Thus, comparing with \equ{wwwdef} we have $\|g\|_{\WWW}=\|g\|_\real$.)
Assume $g_\l\in \WWW_{[-\R ,0]}$ and $g_\r \in \WWW_{[0,\R ]}$ and let
$$
\eqalign{
\E\,=\,
\bigl\{&u\in \CC^2([-\R ,\R ])~:~
|u''|_{[-\R ,\R ]}\le G~,\cr
&\|u-g_\l\|_{[-\R ,0]}\le\epsilon ~,
\|u-g_\r \|_{[0,\R ]}\le\epsilon 
\bigr \}~.\cr
}
$$
\CLAIM{Theorem}{connect}There are a $K$ (depending only on $\epsilon
$ and $G$), and functions
$g_1,\dots,g_N\in\WWW_{[-\R ,\R ]}$ satisfying
$$
\eqalign{
g_i(-\R )\,=\,g_\l(-\R )~,&\quad
g_i(\R )\,=\,g_\r (\R )~,
}
\EQ{smooth}
$$
for $i=1,\dots,N$,
such that the following
holds:
For every
$u\in \E$,
there is a $j=j(u)\in\{1,\dots,N\}$ such that
$$
\|u-g_j\|_{[-\R ,\R ]}\,\le\,\epsilon ~.
$$

\LIKEREMARK{Definition}We say a set 
$\{g_1,\dots,g_K\}$ of functions $g_i\in\WWW$ {\em$\epsilon $-covers a
set $\FF$ of 
$\WWW$ functions on $I$} if for every $g\in \FF$ there is a
$k\in\{1,\dots,K\} $ for which
$$
|g-g_k|_I\,\le\,\epsilon ~,\quad{\rm and}\quad
|g'-g'_k|_I\,\le\,\epsilon ~.
$$
\CLAIM{Corollary}{both}Assume that a collection $\FF$ of $\CC^2$ functions is
given on $[-L,L']$ and assume that each $f\in\FF$ satisfies
$|f|_{[-L,L']}\le \alpha $, $|f'|_{[-L,L']}\le \beta $, and
$|f''|_{[-L,L']}\le\gamma$. There are constants $R$, $\epsilon _0>0$
and a family of 
constants $K_\epsilon$ (depending only on $\alpha $, $\beta $, and
$\gamma$)
such that the following holds for any $L,L'>R$ and any $\epsilon \le
\epsilon _0$:
If $\FF|_{[-L,0]}$ and
$\FF|_{[0,L']}$ can be $\epsilon $-covered by $S$, (resp. $S'$) functions
in $\WWW_{[-L,0]}$ and $\WWW_{[0,L']}$
respectively, then $\FF|_{[-L,L']}$ can be $\epsilon $-covered by no more than
$S\cdot S'\cdot K_\epsilon $  functions in $\WWW_{[-L,L']}$.

\LIKEREMARK{Proof of \clm{connect}}We will first find finite
constants $a$, $b$, $c$ ($>1$)
with the following property:
Fix $g_\l$ and $g_\r $ and assume
$E\equiv \E\ne\emptyset$ (that is, there is a connecting
function in an $\epsilon $-neighborhood of $g_\l$ and $g_\r$). 
We claim one can construct a $W^{2,\infty } $ function $g$ for which
the following 
inequalities hold:
$$
\|g-g_\l\|_{[-\R ,0]}\,\le\,a\epsilon ~,\quad
\|g-g_\r \|_{[0,\R ]}\,\le\,b\epsilon ~,\quad
|g''|_{[-\R ,\R ]}\,\le\,c+G~.
\EQ{lemme1}
$$
Furthermore, $g$ will satisfy
$$
g(-\R )\,=\,g_\l(-\R )~,
\quad g(\R )\,=\,g_\r (\R )~.
\EQ{boundary}
$$

In other words, this is in principle a good approximation, which in addition
matches {\em exactly} at the boundary, but 
the bound has deteriorated to $a\epsilon $ and $b\epsilon $
and $a$ and $b$ might be larger
than 1. The point of
\clm{connect} and \clm{both} is that $a$ (and $b$) can be
pushed down to $1$ by increasing the number of connecting functions
to a number of functions which {\em does not} depend on
$g_\l$ and $g_\r$.

Fix an arbitrary function $u_0\in E$. We construct a function $g$
which interpolates between $g_\l$ and $g_\r $, using $u_0$ as a bridge.
Let $\psi$ be a
$\CC^\infty $ function, $0\le\psi(x)\le 1$ satisfying $\psi (x)=0$ for
$x<R-3$ and $\psi (x)=1$ for $x\ge R$.
We define $g$ by
$$
g(x)\,=\, 
u_0(x) -\psi (x)\cdot\bigl(u_0(\R )-g_\r (\R )\bigr)
-\psi (-x)\cdot\bigl(u_0(-\R )-g_\l(-\R )\bigr)~.
\EQ{gdef}
$$
This function is clearly continuously differentiable since
$u_0$ is continuously differentiable.
Let $I=[0,\R ]$. From \equ{gdef} we find for $x\in I$:
$$
g(x)-g_\r (x)\,=\,\bigl(u_0(x)-g_\r (x)\bigr )-\psi (x)\cdot\bigl(u_0(\R )-g_\r (\R )\bigr)~,
$$
and therefore,
$$
\eqalign{
|g-g_\r |_I\,&\le\,|u_0-g_\r |_I +|u_0-g_\r |_I\cdot|\psi |_I\,\le\,2\epsilon ~,\cr
|g'-g'_\r |_I\,&\le\,|u_0'-g'_\r |_I+|u_0-g_\r |_I\cdot |\psi '|_I
\,\le\,\epsilon \bigl(1+|\psi' |_I\bigr )~.\cr
}
$$
The negative $x$ are handled in the same way.
Finally, the last inequality of Eq.\equ{lemme1} follows
at once from Eq.\equ{gdef}.
We note that by the construction in Eq.\equ{gdef}, the boundary
condition
\equ{boundary} is fulfilled.

\LIKEREMARK{Definition}We denote by $\FABC$ the set of $\CC^2$
functions defined by
$$
\eqalign{
\FABC\,&=\,\{ f ~:~ 
|f(\pm \R )|\le\epsilon ~,
|f|_{[-R,R]} \le A~, |f'|_{[-R,R]} \le B~, |f''|_{[-R,R]}\le C\}~.
}$$
We shall need later the sets
$$
\FABCnull\,=\,\{ f ~:~ |f(0)|\le\delta  ~,~~ 
|f|_{[0,\xi]} \le A~, |f'|_{[0,\xi]} \le B~, |f''|_{[0,\xi]}\le C\}~.
$$

Let $u\in \E$ and define $g$ as in Eq.\equ{gdef}. 
If we let 
$f=u-g$, then by
Eq.\equ{lemme1}, and the construction of
$g$, we see that $f\in \FABC$, with $A=B=(a+b)\epsilon $ and
$C=c+G$.

We will now use the following bound on $\FABC$.
\CLAIM{Proposition}{approx}Fix $0<\epsilon\le A$, $B\ge 0$, and $C\ge 1$.
There is
a {\bf finite} set $H$ of $\WWW $ funct\-ions $H=\{h_1,\dots,h_N\}$,
which $\epsilon $-covers $\FABC$ on $[-\R ,\R ]$ and which furthermore satisfies
$$
h_j(\pm\R )\,=\,0~,
\EQ{rand}
$$
for $j=1,\dots,N$.

Using \clm{approx} we can complete the proof of \clm{connect}.
Given $g_\l$ and $g_\r $ as above, we construct first a $g$ as in
Eq.\equ{gdef}. When $u\in E$, then $f=u-g$ 
is in $\FABC$ by the bounds Eq.\equ{lemme1} and the equality \equ{boundary}.
Thus, by \clm{approx} the $f$ are $\epsilon$-covered
by the $N$ functions $\{h_1,\dots,h_N\}$. 
Define now $u_i=h_i+g$, $i=1,\dots,N$,
and then the set $\E$ of functions $u$ is $\epsilon$-covered by the
$u_i$, since
$$
u-u_i \,=\, (u-g) - (u_i-g) \,=\, f-h_i~,
$$
and we have just stated that the $f$ are $\epsilon $-covered by a
finite number of $u_i$.
Furthermore, the $h_i$ vanish at the boundary of $[-\R ,\R ]$.
Thus, we have interpolated between $g_\l$
and $g_\r $, with $N$ functions in $\WWW$ which $\epsilon
$-cover the original set. The proof of \clm{connect} is complete.\QED
The
corollary then follows at once since the factor $N$ does not depend on
the choice of $g_\l$ and $g_\r $ (except that the bound is too
pessimistic in case $\E$ happens to be empty).

\REMARK The difficulty in proving \clm{approx}
lies in the fact that the $h_j$ vanish at the
endpoints while
the functions $f$ in $\FABC$ may be as large as $\epsilon $ near
the boundary,
$|f(\pm\R )|=\epsilon $, so
there is no space near $\pm R$ with which just to construct an open
cover.

The main ingredient to the proof of \clm{approx} is the following
local lemma. Before we formulate it, we assume, without loss of
generality, that $C>1 $. Since we are interested in small
$\epsilon $, we shall also assume $\epsilon <1$.

We introduce two fundamental scales $\xi$ and $\tau $ in our analysis:
$$
\xi \,=\, {\epsilon \over \X  C}~,\quad{\rm and}\quad
\tau \,=\,{\epsilon \over \T  }~.
$$
We will first
consider a (small) interval $J$ whose left endpoint is the origin.
\LIKEREMARK{Definition of $R$}We can now fix $R$ by setting it to
$$
\the\TRplusvalue\,\ge\,R\,\ge\,\TR~, \quad R\,=\, m_* \xi~,
\EQ{r}
$$
where $ m_* $ is an integer. This choice is only good for $\epsilon
\le \X C$ and we leave the trivial modifications for arbitrary $\epsilon $
to the reader.

\CLAIM{Lemma}{onestep}Let $J=[0,\xi]$. There is a
finite set of linear functions of the form $g_j(x)=j\tau x$,
which $\epsilon $-covers $\FABCnull$, for every $\delta
\in[0,\epsilon ]$. One has in fact better bounds: there is
for every $f\in\FABCnull$ a $j$ with $|j|\le {B\over \tau }+2$ for which
$$
|f-g_j|_J \,\le\,\max(\delta,\epsilon^2\nu)  ~,\quad {\sl and}\quad
|f'-g'_j|_J\,\le\,\epsilon{3\over \T } ~,
$$
and furthermore, at the right endpoint, one has
$$
|f(\xi)-g_j(\xi)| \,\le\,\max\bigl (\delta -\mu\epsilon
^2,\epsilon^2\nu\bigr )~,
\EQ{rightend}
$$
where 
$$
\nu\,=\,{1\over \TNU C}~,\quad \mu\,=\,{1\over \TMU C}~.
\EQ{numu}
$$

\LIKEREMARK{Proof of \clm{onestep}}This is just a construction of the
``right'' $j$, followed by some verifications.
Note first that if $f\in\FABCnull$, then we have
$$
\eqalign{
f(x)\,&=\, f(0) + x f'(0) +x^2 v(x)~,\cr
f'(x)\,&=\,f'(0) + x w(x)~,\cr
}
$$
with $|v|_{[0,\xi]}\le C/2$ and $|w|_{[0,\xi]}\le C$.
We will pursue the proof for the case when $f(0)\ge0$, the other case
is handled by symmetry. We define
$$
j\,=\, \left [ {f'(0)\over \tau  }+2\right ]\,=\, { f'(0)\over
\tau  }+1+ \rho ~,
\EQ{jdef}
$$
with $ \rho\in (0,1]$. Here $[x]= \inf _{n\in\integer, n\ge x}n $ is
the integer part of $x$. Now set $g(x) = c x$, with $c=j\tau $:
$$
c\,=\,f'(0)+\tau +\tau  \rho~.
$$ 
Clearly, $g$ equals one of the $g_j$ of \clm{onestep} if we
take the finite set of $j$ to contain $|j|\le {B\over \tau }+2$.
Next, we estimate the quality of the approximation.
First we have
$$
f'(x)-g'(x)\,=\,f'(x)-c\,=\,f'(0)+xw(x)-f'(0)-\tau -\tau  \rho~.
$$
This leads, for $x\in[0,\xi]$, to
$$
\eqalign{
f'(x)-c \,&\le\, \hphantom{-}C\xi -\tau \,\le\,C\xi \,\le\,{\epsilon
\over \X }~,\cr 
f'(x)-c\,&\ge\,  -C\xi -2\tau \,\ge\, -{\epsilon
\over \X }-{2\epsilon \over \T } ~.\cr
}
$$
We conclude that
$$
|f'-g'_j|_J\,\le\, {3\over \T }\epsilon ~.
$$
We consider next
$f(x)-cx$. We find
$$
f(x)-cx\,=\,f(0)+xf'(0)+x^2 v(x) -xf'(0)-\tau x- \rho
\tau x~.
$$
This leads to the bounds
$$
\eqalignno{
f(x)-cx\,&\le\, f(0) +{C\over 2}x^2-\tau x\,\le\, \delta
+{C\over 2}x^2-\tau x\,=\,
\delta -{\epsilon \over \T }x(1-x{C \over 2 }\cdot {\T\over \epsilon })
~,\NR{fg}
f(x)-cx\,&\ge\, f(0) -{C\over 2}x^2-2\tau x\,\ge\, 
-{C\over 2}x^2-2\tau x~.\NR{fg2}
}
$$
Since we consider only $x\in[0,\xi]$,
we find that $1-x{\T C \over 2\epsilon }\ge1- {\epsilon \over \X C}{\T C \over 2\epsilon }= {1\over 2}$ and therefore
$$
f(x)-cx\,\le\,\delta  ~.
\EQ{upper}
$$
Recall that we deal with the case $f(0)\ge0$. 
Thus, we also get for $x\in
[0,\xi]$,
$$
f(x)-cx\,\ge\,-{C\over 2}{\epsilon ^2\over \X ^2C^2}-{2\epsilon ^2\over
\TX C}\,=\,
-{\epsilon^2 \over 200 C}-{4\epsilon^2 \over 200C }\,=\,-{\epsilon
^2\over \TNU C}\,\equiv\,-\epsilon^2\nu~.
\EQ{lower}
$$
Thus, we conclude that $|f-g_j|_J\le\epsilon  $, provided
$\epsilon\le\TNU C$. 
(This un-intuitive bound comes from having chosen
$\xi=\epsilon /(\X C)$ which is unreasonable when $C\ll 1$.)

We next show that the bound on $f(x)-cx$ is tighter than what we got
so far
when $x=\xi$.
Indeed, we get in this case from Eq.\equ{fg},
$$
f(\xi)-c\xi \,\le\,\delta  -{1\over 2}{\epsilon^2 \over \TX C}\,\equiv\,\delta
-\mu\epsilon ^2~.
\EQ{fg3}
$$
The
assertion Eq.\equ{rightend} follows by combining Eq.\equ{fg3} with
Eq.\equ{lower}.

It remains to see that the set of possible $j$ is finite. Considering
Eq.\equ{jdef} and the fact that $f\in\FABCnull$ we see that $j$ can take
at most $2(B/\tau+2  )+1$ possible values.
The proof of \clm{onestep} is complete.\QED

\LIKEREMARK{Proof of \clm{approx}}This proof is a repeated application
of \clm{onestep}. We retain the assumptions and notations from that
proof.
Let $\eta=\xi\tau={\epsilon ^2\over \TX C}$.
We consider the grid (in the $(x,y)$-plane):
$$
\bigl \{\bigl ( m\xi,n\eta\bigr )~:~  m=-m_*,-m_*+1,\dots, m_*~;~
n=-n_*,\dots,n_*\bigr \} ~,
$$
where
$$
m_*\,=\,\R /\xi~,\quad n_*\,=\, [A/\eta]+1~,
$$
recalling that $R/\xi$ is an integer.
In other words, we cover the range of possible arguments (in $[-\R ,\R ]$)
and values (in $[-A,A]$) of 
$f\in\FABC$ by a fine grid.
Consider now the set of all continuous, piecewise linear functions
$h(x)$, connecting linearly successive lattice points
$(m\xi ,n \eta)$ with $((m+1)\xi,n'\eta)$, with $-m_*\le m< m_*$, $|n|\le
n_*$ and $|n'|\le 
n_*$. Furthermore, we require that $h(-\R )=h(\R )=0$.
There are a finite number
of such functions, namely at most $(2n_*+1)^{2m_*-1}$.

Note that $\eta$ has been chosen in such a way that the slopes of the
straight pieces of $h$ are integer multiples of $\tau $.
We show next that every $f\in\FABC$ is, together with its derivative,
$\epsilon$-close to one of the $h$. 

We begin by constructing the piecewise linear approximation $h$.
We start at the point $x=-\R $, $y=0$, and shift the
origin to this point by defining:
$$
f_0(x) \,=\, f(x+\R)~.
$$
Then $f_0$ is in $F^0_{\epsilon
,A,B,C,2R}\supset \FABCnull$, with $\delta =\epsilon $, and 
by \clm{onestep}, $f_0$ is
approximated by one of the linear functions, say $n_0\tau x$, with 
$n_0=[2+f_0'(0)/\tau ]$ on the interval
$[0,{\xi}]$ (when $f(0)>0$). Note that we also have (when $\epsilon $
is small)
$|f_0(\xi)-n_0\tau \xi|\le \epsilon
-\mu\epsilon ^2$.
We define
$$
h(x)\,=\,n_0\tau (x+R)~, \text{ for } x\in [-R,-R+\xi]~.
$$
Next shift the origin of the
$(x,y)$-plane to $(-\R +\xi,n_0\eta)=\bigl (-\R+\xi,h(\xi)\bigr ),$ and define 
$$
f_1(x) \,=\, f(x+\R -\xi) -h(x+\R-\xi)\,=\,f(x+\R-\xi)-n_0\eta~.
$$
The definition of the first segment of $h$ and the
bounds on $f_0$ 
show that 
$$
f_1\in F^0_{\epsilon-\mu\epsilon ^2 ,A+|n_0|\eta,B,C,2R-\xi}~.
$$

We now apply \clm{onestep} to $f_1$. Note that $f_1$ is not in $\FABCnull$
but in a space with a worse bound on the absolute value. However, the
value of $A$ does {\em not} enter the construction of the proof of
\clm{onestep} and hence is irrelevant for our inductive construction
of $h$. Applying \clm{onestep} to $f_1$, we find the second linear
piece of the function $h$, and get a piecewise linear, continuous
approximation of $f$ on
$[-\R ,-\R +2\xi]$.
The final point of the approximation by $h$ is now $(-\R +2\xi,
n_1\eta)$, and we construct 
$f_2$ by translating the origin to that point. 
Assuming that $\epsilon -2\mu\epsilon ^2>\nu\epsilon ^2$, we conclude
that
$$
f_2\in F^0_{\epsilon  -2\mu\epsilon ^2,A+|n_0+n_1|\eta,B,C,2R-2\xi}~.
$$
Note that the
construction can not ``drift away'' in the $y$-direction, since we
assumed from the outset that $|f|_{[-R,R]}\le A$, and hence the
$y$-translates never exceed $A$ by more than $\epsilon $ (since $h$
{\em is} an approximation to $f$). We continue the construction in the
same way as before, until $x=0$ is reached.
At this point we have achieved the following:
The original function is approximated by the piecewise linear function
$h$ on $J=[-\R,0]$ with the bound
$$
|f-h|_J\,\le\,\epsilon ~,\quad
|f'-h'|_J\,\le\,\epsilon ~.
$$
Furthermore, at the point $x=\xi$ the approximation is really ``good:''
Consider the definition \equ{r} of $R$.
The number of steps from $-R$ to $0$ is $ m_* \ge {\TR
}\cdot{1\over \xi}-1={\TR \cdot\X C\over \epsilon }-1$
and in each step we gain a constant $\mu\epsilon ^2$, as long as $\delta
>\epsilon^2\nu$. Therefore,
$$
|f(0)-h(0)|\,\le\, \max\bigl (\epsilon^2\nu,\epsilon
- m_*\epsilon ^2\mu\bigr )\,=\,\epsilon^2\nu~,
\EQ{good1}
$$
where the last equality follows from
\Ttemp=\TR\multiply\Ttemp by \X
$$
\epsilon ^2\mu m_* \,\ge\,\epsilon ^2\mu( m_*-1) \,\ge\, \epsilon
^2\mu
\bigl ( {\the\Ttemp C\over \epsilon }-2
\bigr )\,=\, 
\epsilon ^2{1\over \TMU C}\cdot\bigl ({\the\Ttemp C\over \epsilon }-2\bigr
)\,\ge\,
2\epsilon  -{2\epsilon ^2\over \TMU C}\,\ge\,\epsilon ~,
\EQ{rho}
$$
when $\epsilon\le 100C$.

We repeat the same construction from the right endpoint, (with
$ m_*-1$ steps, which is also covered by \equ{rho})
obtaining the
piecewise linear function $h$ on the set $J=[\xi,\R]$, and again a
bound, using \equ{numu}:
$$
|f(\xi)-h(\xi)| \,\le\,\max\bigl (\epsilon ^2\nu,\epsilon
-( m_* -1)\epsilon ^2\mu\bigr )\,\le\,\epsilon
^2\nu~.
\EQ{good2}
$$
We complete the definition of $h$ by connecting $\bigl (0,h(0)\bigr )$
linearly with $\bigl (\xi,h(\xi)\bigr )$. Note that it is necessarily a
line segment connecting two of the grid points and so $h$ is one of
the functions we counted earlier.
We need to verify the bounds on $J=[0,\xi]$. It is here that the
Eqs.\equ{good1} and \equ{good2} are relevant. We write
$$
f(x)\,=\,f(0)\cdot(1-{x\over \xi})+f(\xi){x\over \xi}+r(x)~,
$$
and then by the bounds on the second derivative of $f$ we get $
|r|_J\le C\xi^2/8$, and $|r'|_J\le C\xi/2$. Since
$$
h(x)\,=\,h(0)\cdot(1-{x\over \xi})+h(\xi){x\over \xi}~,
$$
we find for $\epsilon\le 800C/21 $,
$$
\displaylines{
|f-h|_J\,\le\,\epsilon ^2\nu + {C\over 8}{\epsilon ^2\over
C^2\X ^2}\,=\,\epsilon ({\epsilon \over \TNU C}+{\epsilon \over 8\cdot \X^2
C})\,\le\,\epsilon ~,\cr
|f'-h'|_J\,\le\,{2\over \xi} \epsilon ^2\nu   + {C\xi\over 2}\,=\,
{2\cdot\X C\over \epsilon }\cdot \epsilon ^2 \cdot{1\over \TNU C}+{C\epsilon
\over 2\cdot\X C}\,=\,{11\over 20}\epsilon \,\le\,\epsilon~.
}
$$
Thus, we have shown the required bound on all of $[-R,R]$.
The piecewise linear,
continuous function obtained in this way will be called $h_f(x)$. It
is clearly one of the functions we constructed. It approximates $f$
and $f'$ on all of $[-\R ,\R ]$. We have thus found a finite family of
piecewise linear functions which $\epsilon $-covers $\FABC$.
The proof of \clm{approx} is complete.\QED

\SUBSECT{ee}{The $\epsilon $-entropy of Kolmogorov and Tikhomirov}

We proceed
as in [CE1], but with a change of topology as explained above.
We have defined in \sec{bounds} the minimum number
$N_L(\epsilon )$ of balls in the norm $\|\cdot\|_{\delta,L,2}$ needed to cover
the attracting set. We also showed in \clm{balls} that
$$
N_{L-A/\epsilon }(\epsilon )\,\le\,\CX{ball}^L N_L(2\epsilon
)~,
\EQ{balls2}
$$ 
with some constants $A$ and $\CX{ball}$ depending only on the
coefficients of the problem \equ{problem}. 
If we iterate Eq.\equ{balls2} $m$ times, we get
$$
N_L(\epsilon  )\,\le\,\CX{ball}^{L+A/\epsilon }\CX{ball}^{L+2A/\epsilon }\cdots
\CX{ball}^{L+mA/\epsilon } N_{L+A/\epsilon +2A/\epsilon
+\cdots+mA/\epsilon }(2^m\epsilon )~.
\EQ{balls5} 
$$
In \equ{good0} we have shown that
there is a constant
$\CX{77}'$ which bounds the radius of $\GG$ in $\HH_{\delta,\loc,2}$.
(The bound $\CX{77}$ in \equ{good0} was for $\HH_{\alpha ,\loc,2}$.)
Therefore, {\em one} ball of
radius $\CX{77}'$ suffices to cover $\GG|_{[-L,L]}$. Choosing
$m=m(\epsilon )$ in such a way that $2^m\epsilon >\CX{77}'$, we conclude
that $\GG|_{[-L,L]}$ can be covered by a {\em finite} number of
balls in $\HH_{\delta ,\loc,2}$.\footnote{${}^1$}{The argument used
here is more elegant than the one used in [CE1]. We thank
Y. Colin de Verdi\`ere for suggesting it.}
\REMARK This argument does not use compactness of $\GG_{[-L,L]}$ and
does not prove it either. It is here that our method differs from that
of Feireisl[F]. He shows that the intersection of the
$\Phi^T(\GG|_{[-L,L]})$ is compact, whereas our approach shows that
$\GG|_{[-L,L]}$ itself can be covered by a finite number of balls.

We define similarly
for any interval $I$, the
minimal number $M_{I}(\epsilon )$ of
balls needed to cover $\GG$ in the topology
$\|\cdot\|_{\WWW_I }$. By the Sobolev inequality from
\clm{Sobolev}, we see that 
$$
\|u\|_{\WWW_{[-L,L]}}\,\le\,\CX{Sobolev} \|u\|_{\delta,L,2}~.
$$
Therefore, if $\GG$ can be covered by $N_L(\epsilon
/\CX{Sobolev})$ balls of radius $\epsilon /\CX{Sobolev}$ in the norm
$\|\cdot\|_{\delta,L,2}$, it can obviously be covered by the same
number of balls of radius $\epsilon $ in the norm $\|\cdot\|_{\WWW_{[-L,L]}}$.
Thus we have
$$
M_{[-L,L]}(\epsilon )\,\le\,N_L(\epsilon /\CX{Sobolev})~.
\EQ{compare} 
$$
We now apply \clm{both}. We first note that by Eqs.\equ{good3} and
\equ{good4} it is adequate to consider functions with {\em bounded}
second derivative. (In fact this is the only place where these higher
derivatives are needed.) Thus, we can apply   
\clm{both} and we conclude that for two
intervals $I_1$ and $I_2$, one has
$$
M_{I_1\cup I_2}(\epsilon )\,\le\, M_{I_1}(\epsilon )M_{I_2}(\epsilon
)K_\epsilon ~.
\EQ{sub} 
$$
Thus, we have established submultiplicativity (in $I$) and finiteness of
$M_{I}(\epsilon )$. Furthermore, from the construction of
$m$ in \equ{balls5} with $2^m\epsilon >\CX{77}' $, we find by
choosing the minimal such $m$:
$$
M_I(\epsilon
)\,\le\,\CX{ball}^{\C{6}\log(\epsilon^{-1})|I|+\log(\epsilon
^{-1})\epsilon ^{-1} \C{66}A}~.
$$ 
Using this bound and \equ{sub}, we get convergence and a bound on the $\epsilon
$-entropy $H_\epsilon (\GG)$: 
\CLAIM{Theorem}{ee}The $\epsilon $-entropy 
 per unit length of $\GG$ in
$\WWW$ exists and is bounded by
$$
H_\epsilon (\GG)\,=\,\lim_{L\to\infty }{1\over L}\log \bigl(M_
{[-L,L]}(\epsilon )\bigr )\,\le\,\C{ee}\log(1/\epsilon  )~.
\EQ{ee}
$$

\SUBSECT{top}{Existence of the topological entropy per unit length}

This material is taken from [CE2], and we introduce it without proofs,
just to show what follows from the bounds of the previous sections.

For any $\epsilon>0$ and any interval $I$ in $\real$,
we define $\wqeps$ as the set of all finite covers of $\attra$
by open sets in $\WWW_I$
of diameter at most $\epsilon$. 
Note that by the argument of \sec{ee}, such a finite cover exists.
Note also that elements of $\WWW_I$ are pairs of functions $(u,v)$ and
that the topology is $W^{1,\infty }_{I}$ on the $u$-component and
$\L^\infty (I)$ of the $v$-component.

Let $\tau>0$ be a fixed time step, and let $T=n\tau$ with
$n\in{\integer}$.
\LIKEREMARK{Definitions}Let ${\cal U}\in \wqeps$. 
We say that two elements $A_{1}$ and
$A_{2}$  in $\attra$ are $(\calu ,T)$-{\em separated}
if there is at least one $k\in\{0,\dots,n\}$ for which the 
points $\Phi^{k\tau}(A_{1})$
and  $\Phi^{k\tau}(A_{2})$ do not belong to the same atom of $\cal U$.
We define
$$
N_{T,\tau}({\cal U })
$$
to be the largest number of elements which are pairwise $(\calu
,T)$-separated
(and considered with time-step $\tau $.) Note that
this number is finite since it is at
most $({\rm Card}~ {\cal U})^{2T/\tau}$.
Finally, we define
$$
N_{I,T,\tau,\epsilon}=\inf_{{\cal U}\in \wqeps}
N_{T,\tau}({\cal U })~.
$$
\CLAIM{Lemma}{toto}(Lemma 2.1. of [CE2]). Let $I_{1}$ and $I_{2}$ be
two disjoint
intervals (perhaps with common boundary) and let $I=I_1\cup I_2$.
The functions  $N_{I,T,\tau,\epsilon}$ satisfy the
following bounds: There is a constant $C=C(\epsilon )$ such that:
\item{i)} $N_{I,T,\tau,\epsilon}$
is non-increasing in $\epsilon$.
\item{ii)} $N_{I,T_{1}+T_{2},\tau,\epsilon}
\le N_{I,T_{1},\tau,\epsilon}~N_{I,T_{2},\tau,\epsilon}$.
\item{iii)} $N_{I_{1}\cup I_{2},T,\tau,\epsilon}
\le C\,N_{I_{1},T,\tau,\epsilon}~N_{I_{2},T,\tau,\epsilon}$.

\REMARK It is important here that $C(\epsilon)$ does {\em not} depend
on the lengths of $I_1$ and $I_2$.
\PROOF The properties i) and ii) are shown exactly as in
[CE2]. However, the proof and the statement of iii) are now modified
since we consider the topology of $\WWW$.

In order to prove iii), we consider ${\cal U}_{1}\in {\cal
W}_{I_{1}}^{\epsilon}$ and ${\cal U}_{2}\in {\cal W}_{I_{2}}^{\epsilon}$.
Since we are using the $\WWW $ norm we have
$$
N_{T,\tau}({\cal U}_{1}\cap{\cal U}_{2})
\,\le\, N_{T,\tau}({\cal U}_{1})
~N_{T,\tau}({\cal U}_{2})~.
$$
We also have easily
$$
{\cal W}_{I_{1}}^{\epsilon}\cap {\cal W}_{I_{2}}^{\epsilon}
\subset {\cal W}_{I_{1}\cup I_{2}}^{\epsilon}~.
$$
The claim iii) now follows easily.\QED

\LIKEREMARK{Remark}Henceforth, we shall work with domains which are
intervals $I_L=[-L,L]$.
\CLAIM{Theorem}{exist}The following limit exists
$$
h\,=\,\lim_{\epsilon\to0}\lim_{L\to\infty}{1\over L}
\lim_{T\to\infty}{1\over T}
\log N_{I_{L},T,\tau,\epsilon}~.
\EQ{lim}
$$
Moreover, $h$ does not depend on $\tau$. It is called the topological
entropy per unit volume of the system.

\PROOF The proof is given in [CE2].
\REMARK It also follows from \sec{ee} that $h$ is bounded.
\SUBSECTION{Sampling}

The results we describe in this section are, on the surface, the same
as those obtained in [CE2]. This means that by {\em discrete} sampling
of the signal in
a space-time region 
$$
[-L-A\log(1/\epsilon ),
L+A\log(1/\epsilon )]\times[0,\tau _* \log(1/\epsilon )]~,\EQ{region} $$ 
the function observed can be determined to precision $\epsilon $ {\em
everywhere} on the
interval $[-L,L]$ at time $\tau _*\log(1/\epsilon )$.

In the current context this result can be worked out in detail in the
following sense:
Assume that two solutions $u_1$ and $u_2$ and their first and second
space derivatives (as well as $\partial _t u_1$ and $\partial _t u_2$ and
their first derivatives) coincide to within $\epsilon $ in the region \equ{region}
on a space-time grid with
mesh $\OO(1/k_*)\times \OO(\tau _*)$. Then one can conclude that
$$\eqalign{\|u_1(&\tau _*\log(1/\epsilon),\cdot
)-u_2(\tau_*\log(1/\epsilon),\cdot 
)\|_{\WWW_{[-L,L]}}\cr+\|\partial _tu_1(&\tau _*\log(1/\epsilon),\cdot
)-\partial _t u_2(\tau_*\log(1/\epsilon),\cdot
)\|_{\L^\infty _{[-L,L]}}\,\le\, \C{last}\epsilon~,}
$$
for some universal constant $\CX{last}$. This result allows, in
principle, to reconstruct the $K_2$-entropy. 

In our view, the result sketched above is somewhat unsatisfactory, and
its clarification needs further work. Namely, we would like to be able
to make positive statements based on sampling {\em only} function
values, and not their derivatives, in particular, not the second
derivative. (They are needed to bound the difference in $\WWW$.)
Indeed, a quick inspection of properties of the Bernstein class shows
that we have no reasonable bound on $\SS_L( \QQ_{k_*}f) - \SS_L( f)$ in
$\WWW$ if we
only have information about the function and not about its derivatives. 
\SUBSECTIONNONR{Acknowledgments}

We have profited from useful
discussions with Th.~Gallay and J. Rougemont. This work was partially
supported by the Fonds National Suisse.
\SUBSECTIONNONR{References}

\ref
\no{B}
\by{Boas, R.P.}
\book{Entire Functions}
\publisher{New York, Academic Press}
\yr{1954}
\endref
\ref
\no{CE1}
\by{Collet, P. and Eckmann, J.-P.}
\paper{Extensive properties of extended systems}
\jour{Commun.Math.Phys.}
\vol{200}
\pages{699--722}
\yr{1999}
\endref

\ref
\no{CE2}
\by{Collet, P. and Eckmann, J.-P.}
\paper{Topological
entropy per unit volume in parabolic PDE's}
\jour{Nonlinearity}
\vol{12}
\pages{451--473}
\yr{1999}
\endref

\ref
\no{F}
\by{Feireisl, E.}
\paper{Bounded, locally compact global attractors for semilinear damped wave
equations on $\real^ N$} 
\jour{Differ. Integral Equ.} 
\vol{9}
\pages{1147--1156}
\yr{1996}
\endref
\ref
\no{KT}
 \by{Kolmogorov A.N. and Tikhomirov, V.M.}
 \paper{$\epsilon$-entropy and $\epsilon $-capacity of sets in
 functional spaces\footnote{${}^1$}{The version in this collection is
 more complete than the original paper of Uspekhi Mat. Nauk, {\bf 14},
 3--86 (1959).}}
 \inbook{Selected Works of A.N. Kolmogorov, Vol III}
 \bybook{Shirayayev, A.N., ed.}
 \publisher{Dordrecht, Kluver}
 \yr{1993}
\endref
\ref
\no{M1}
\by{Mielke, A.}
\paper{The complex Ginzburg-Landau equation on large and unbounded
domains: sharper bounds and attractors}
\jour{Nonlinearity}
\vol{10}
\pages{199--222}
\yr{1997}
\endref

\ref
\no{M2}
\by{Mielke, A.}
\paper{private communication}
\endref

\ref
\no{MS}
\by{Mielke, A. and Schneider, G.}
\paper{Attractors for modulation equations on unbounded domains --
existence and 
comparison}
\jour{Nonlinearity}
\vol{8}
\pages{743--768}
\yr{1995}
\endref

\ref
\no{R}
\by{Ruelle, D.}
\book{Statistical Mechanics}
\publisher{New York: Benjamin}
\yr{1963}
\endref
\bye